\documentclass[12pt]{amsart} 
\usepackage[dvips]{graphics}
 \usepackage{amssymb, amsmath, mathrsfs}
  \input xy
  \xyoption{all}
  \DeclareMathAlphabet{\mathpzc}{OT1}{pzc}{m}{it}

 \newtheorem{theorem}{Theorem}[section]
 \newtheorem{proposition}[theorem]{Proposition}
 \newtheorem{corollary}[theorem]{Corollary}
 \newtheorem{lemma}[theorem]{Lemma}

 \newtheorem{conjecture}[theorem]{Conjecture}
 \theoremstyle{definition}
 \newtheorem{definition}[theorem]{Definition}

 \theoremstyle{remark}
 \newtheorem{remark}[theorem]{Remark}
 
 \newtheorem{remarks}[theorem]{Remarks}

\newcommand{\CA}{{\mathcal A}}
\newcommand{\CB}{{\mathcal B}}

\newcommand{\CE}{{\mathcal E}}
\newcommand{\CF}{{\mathcal F}}
\newcommand{\CG}{{\mathcal G}}
\newcommand{\CH}{{\mathcal H}}
\newcommand{\CI}{{\mathcal I}}

\newcommand{\CM}{{\mathcal M}}

\newcommand{\CO}{{\mathcal O}}
\newcommand{\CP}{{\mathcal P}}

\newcommand{\CS}{{\mathcal S}}
\newcommand{\CT}{{\mathcal T}}

\newcommand{\CV}{{\mathcal V}}
\newcommand{\CW}{{\mathcal W}}

\newcommand{\CZ}{{\mathcal Z}}

\newcommand{\fh}{{{\mathfrak h}}} 
 
\newcommand{\fp}{{{\mathfrak p}}} 
\newcommand{\fg}{{{\mathfrak g}}} 
\newcommand{\fb}{{{\mathfrak b}}}

\newcommand{\fhd}{\fh^\star}

\newcommand{\SB}{{\mathscr B}}

\newcommand{\SL}{{\mathscr L}}
\newcommand{\SM}{{\mathscr M}}
\newcommand{\SN}{{\mathscr N}}

\newcommand{\DP}{{\mathbb P}}
\newcommand{\DR}{{\mathbb R}}
\newcommand{\DZ}{{\mathbb Z}}

\newcommand{\DN}{{\mathbb N}}

\newcommand{\ch}{{\operatorname{char}}}
\newcommand{\im}{{\operatorname{im}}}

\newcommand{\Sh}{{\mathcal{SH}}}

\newcommand{\End}{{\operatorname{End}}}

\newcommand{\Hom}{{\operatorname{Hom}}}
\newcommand{\id}{{\operatorname{id}}}

\newcommand{\udim}{{\operatorname{\underline{dim}}}}

\newcommand{\Gr}{{\operatorname{Gr}}}

\newcommand{\supp}{{\operatorname{supp}\,}}
\newcommand{\catmod}{{\operatorname{-mod}}}

\newcommand{\bimod}{\text{--mod--}}
\newcommand{\inj}{{\hookrightarrow}}

\newcommand{\newtimes}{\times}

\newcommand{\Nabla}{\bigtriangledown}

\newcommand{\eval}{{\operatorname{ev}}}
\newcommand{\Stab}{{\operatorname{Stab}}}

\newcommand{\Loc}{{\SL}}
\newcommand{\on}{{T_{on}}}
\newcommand{\out}{{T_{out}}}

\newcommand{\ol}{\overline}

\newcommand{\linie}{{\,\text{---\!\!\!---}\,}}
\newcommand{\llinie}{{\text{---\!\!\!---\!\!\!---}}}

\begin{document}

\pagenumbering{arabic}
\title[]{The combinatorics of Coxeter categories}
\author[]{Peter Fiebig}
\address{Mathematisches Institut, Universit{{\"a}}t Freiburg, 79104 Freiburg, Germany}
\email{peter.fiebig@math.uni-freiburg.de}
\begin{abstract} We present an alternative construction of Soergel's category of bimodules associated to a reflection faithful representation of a Coxeter system. We show that its objects can be viewed as sheaves on the associated moment graph. We introduce an exact structure and show that the ``special bimodules'' are the projective objects. Then we construct the indecomposable projectives  by both a global and a local method, discuss a version of the Kazhdan-Lusztig conjecture and prove it for universal Coxeter systems. 
\end{abstract}
\maketitle
\section{Introduction}

Coxeter systems abound in representation theory and geometry and can often be used to describe categorical structures. Conversely, some of the most intriguing combinatorial problems for Coxeter systems can, as yet, only be proven by ``categorification''.  The main objective of this article is to study such a categorification.

Let $V$ be a {\em reflection faithful} representation of a Coxeter system $(\CW,\CS)$ over a field $k$ of characteristic $\neq 2$ (see Definition \ref{refrep}).  Let $S=S_k(V^\ast)$ be the symmetric algebra over $V^\ast=\Hom_k(V,k)$, graded such that $\deg V^\ast=2$. Let $\Lambda\subset V$ be an orbit of $\CW$ and choose a partial order on $\Lambda$ such that two elements are comparable if they differ by a reflection in $\CW$. To such data we assign an additive $S$-category $\CV=\CV(\Lambda)$ together with an exact structure.  

We consider three different realizations of the category $\CV$. The first and maybe the easiest is the following. There is an associative and commutative $S$-algebra $\CZ$ assigned to the action of $\CW$ on the orbit $\Lambda$. The partial order is used to define a cofiltration on each $\CZ$-module in a functorial way. Then $\CV$ is the subcategory of objects such that each quotient in the filtration is graded free over $S$. 

The second realization is similar and is due to Soergel. Suppose that $k$ is infinite and identify $S$ with the space of polynomial functions on $V$. We consider $S$-bimodules supported on a finite union of twisted diagonals in 
$V\newtimes V$ (cf.~ Section \ref{S-bimod}). The length function on $\CW$ induces a support filtration on such bimodules in a functorial way. Let $\CF_\Nabla$ be the category of $S$-bimodules whose subquotients are isomorphic to a direct sum of copies of the polynomial functions on the respective diagonals, shifted in degree. 

Assume that $\Lambda\subset V$ is a regular orbit, i.e.~ that the stabilizer of each point is trivial, and choose $v\in\Lambda$. Identify $\CW$ and $\Lambda$ via the action of $\CW$ on $v$. We will show that the category $\CV(\Lambda)$, that is constructed using the partial order induced by the Bruhat order, is equivalent to $\CF_\Nabla$. 

The third realization is different and more local in nature. From the partially ordered orbit $\Lambda$ (not necessarily regular) we construct a {\em moment graph} $\CG$, i.e.~ a directed graph without cycles and double edges and a labelling of each edge by a one-dimensional subspace of $V^\ast$. In \cite{BMP} Braden and MacPherson defined the notion of a sheaf on such a graph. We construct $\CV(\CG)$ as the category of all sheaves that are generated by global sections and are flabby and have the property that the sections on each upwardly closed subgraph form a graded free module over $S$.

In \cite{Fie05} it was shown that the first and the third realization lead to equivalent categories. The first result in this article is the equivalence to Soergel's construction. 

The exact structure on $\CV$ is naturally associated to the partial order on $\Lambda$ and was defined in \cite{Fie05} (cf.~ also \cite{Dyer2} for an analogous exact structure).  It allows us to introduce notions common to abelian categories, such as short exact sequences, exact functors, or projective objects. Let $\CP\subset\CV$ be the full subcategory of projective objects. We show that the subcategory $\CB\subset\CF_\Nabla$ of ``special bimodules'', defined in \cite{Soe04}, corresponds to $\CP$ under the equivalence $\CF_\Nabla\cong\CV$. This is done using a well-known construction of projectives by translation functors. We also show that the indecomposable projectives correspond to the {\em Braden-MacPherson sheaves} on $\CG$.

It is conjectured in \cite{Soe04} that the graded characters of the indecomposable objects in $\CB$ are given by the Kazhdan--Lusztig self-dual elements in the Hecke algebra.  We use the equivalence of the global and the local realizations to support the conjecture and to prove it in the case of universal Coxeter systems. I hope that the comparison of the global and local viewpoints turn out to be valuable beyond the results of this article.

\section{A category associated to an orbit of a  Coxeter system}\label{catcox}

\subsection{Reflection representations}
Let $(\CW,\CS)$ be a not necessarily finite Coxeter system and $\CT\subset\CW$ the set of reflections in $\CW$, i.e.~ the orbit of $\CS$ under conjugation. Choose a field $k$ of characteristic $\neq 2$ and let $V$ be a finite dimensional representation of $\CW$ over $k$. For $w\in\CW$ we denote by $V^{w}\subset V$ the set of $w$-fixpoints. The following is Soergel's definition of {\em spiegelungstreu}. 

\begin{definition}[\cite{Soe04}, Definition 1.5]\label{refrep} By a {\em reflection faithful representation} of $(\CW,\CS)$ we mean a faithful, finite dimensional representation $V$ of $\CW$ such that, for each $w\in\CW$, the subspace $V^w$ is a hyperplane in $V$  if and only if $w\in\CT$. 
\end{definition}

The {\em geometric representation} of $(\CW,\CS)$ (cf.~ \cite{Hum}) might not be reflection faithful. In \cite{Soe04} a reflection faithful representation was constructed over the field $\DR$ for each Coxeter system $(\CW,\CS)$. In the crystallographic cases it coincides with the representation of $\CW$ on the Cartan subalgebra of the associated Kac-Moody algebra. If $(\CW,\CS)$ is crystallographic and if we denote by $L$ its root lattice, then the induced representation on $L\otimes_\DZ k$ being reflection faithful depends on the characteristic of $k$. 
 
From now on suppose that $V$ is reflection faithful. For a reflection $t$ let $V^{-t}$ be the $t$-eigenspace for the eigenvalue $-1$.
\begin{lemma} \label{rfrep} Let $s,t\in\CT$. Then
\begin{enumerate}
\item  $V^s=V^t\Leftrightarrow   V^{-s}=V^{-t} \Leftrightarrow s=t$.
\item Let $x\in V$ and suppose $s.x=t.x\neq x$. Then $s=t$.
\end{enumerate}
\end{lemma}
\begin{proof}
Part (1) is proven in \cite{Soe04}. Note that for each $s\in\CT$ and $x\in V$,  $x-s.x\in V^{-s}$, and hence (2) follows from (1) and the fact that $V^{-s}$ is one-dimensional for $s\in\CT$. 
\end{proof}

Let $S=S_k(V^\ast)$ be the symmetric algebra of $V^\ast=\Hom_k(V,k)$, endowed with the unique algebra $\DZ$-grading given by setting $V^\ast$ in degree $2$.  For  a graded $S$-module $M$ and $n\in\DZ$ denote by $M_{\{n\}}$ its homogeneous component of degree $n$. In the following we will only consider $\DZ$-graded $S$-modules, and a morphism $f\colon M\to N$ of graded $S$-modules will always be of degree zero, i.e.~ such that $f(M_{\{n\}})\subset N_{\{n\}}$ for all $n\in\DZ$.

\subsection{The structure algebra}

For each $t\in\CT$ let $\alpha_t\in V^{\ast}$ be a non-trivial linear form vanishing on the hyperplane $V^t$. Then $\alpha_t$ is well-defined up to a non-zero scalar. Moreover, $\alpha_t$ and $\alpha_s$ are linearly independent if $s\neq t$ by Lemma \ref{rfrep}. 

Let $\Lambda\subset V$ be an orbit of $\CW$ and define 
$$
\CZ=\CZ(\Lambda):=\left\{\left. (z_x)\in\prod_{x\in\Lambda}S\, \right|\,\begin{matrix} z_x\equiv z_{t.x}\mod \alpha_t \\ \text{for all $x\in\Lambda$ and $t\in\CT$} \end{matrix}\right\}.
$$
$\CZ$ is a commutative, associative, $\DZ$-graded $S$-algebra. For any subset $\Omega\subset\Lambda$ let $\CZ^\Omega$ be the image of $\CZ$ under the projection $\prod_{x\in\Lambda}S\to \prod_{x\in\Omega}S$ with kernel $\prod_{x\in\Lambda\setminus \Omega}S$. We write $\CZ^{x_1,\dots,x_n}$ instead of $\CZ^{\{x_1,\dots,x_n\}}$. 

\begin{definition}
Let $\CZ\catmod^f$ be the category of graded $\CZ$-modules that are finitely generated over $S$ and torsion free over $S$, and such that the action of $\CZ$ factors over $\CZ^\Omega$ for a {\em finite} subset $\Omega\subset\Lambda$. Let the morphisms be the graded $\CZ$-morphisms of degree zero. 
\end{definition}

Choose $v\in\Lambda$. The action of $\CW$ on $v$ gives an identification $\CW/ \Stab_\CW(v)=\Lambda$ that we fix. For $w\in\CW$ we denote by  $\ol w$ its image in $\CW/ \Stab_\CW(v)$. Let $\CW$ act on $V^\ast$ by $(w.\lambda)(v)=\lambda(w^{-1}.v)$ for $w\in\CW$, $\lambda\in V^\ast$ and $v\in V$. For $\alpha\in V^\ast$ such that $\Stab_\CW(v)\subset \Stab_\CW(\alpha)$ define $\sigma(\alpha)=(\sigma(\alpha)_{\ol w})\in\prod_{\ol w\in\CW/\Stab_\CW(v)} S$ by $\sigma(\alpha)_{\ol w}=w.\alpha$. 
\begin{lemma} 
\begin{enumerate}\label{loccen}
\item We have $\sigma(\alpha)\in\CZ(\Lambda)$. 
\item Choose $x\in\CW$ and $t\in\CT$ such that $\ol{tx}\neq \ol x$ in $\CW/ \Stab_\CW(v)$. Then $\CZ^{\ol x,\ol{tx}}=\left\{(z_{\ol x},z_{\ol{tx}})\in S\oplus S\mid z_{\ol x}\equiv z_{\ol{tx}}\mod \alpha_t\right\}$. 
\end{enumerate}
\end{lemma}
\begin{proof} Let $w,tw\in\CW$ with $\ol w\neq \ol{t w}$ in $\CW/ \Stab_\CW(v)$.  Then $w.\alpha-tw.\alpha$ vanishes on $V^t$, so it must be a multiple of $\alpha_t$. Hence (1). 

We will now prove (2). It is clear that $\CZ^{\ol x,\ol{tx}}$ is contained in the space on the right hand side. It is also clear that $(1,1)\in\CZ^{\ol x,\ol{tx}}$. So we only have to show that $(\alpha_t,0)\in \CZ^{\ol x,\ol{tx}}$. Choose $\alpha\in V^\ast$ with $\Stab_\CW(v)\subset\Stab_\CW(\alpha)$ and $x^{-1}tx\not\in\Stab_\CW(\alpha)$. Then $\sigma(\alpha)_{\ol x}\neq\sigma(\alpha)_{\ol{tx}}$, hence $\sigma(\alpha)-\sigma(\alpha)_{\ol{tx}}\cdot 1$ is an element in $\CZ$ whose image in $\CZ^{\ol x,\ol{tx}}$ is a multiple of $(\alpha_t, 0)$.
\end{proof}

\begin{remark}\label{sigmareg} For a regular orbit $\Lambda$ we get a linear map $\sigma\colon V^\ast\to\CZ$ that extends to an algebra homomorphism $\sigma\colon S\to \CZ$.  
 \end{remark}
\subsection{Quasi-finiteness}
Let $Q$ be the quotient field of $S$. Define $N_Q:=N\otimes_S Q$ for any $S$-module $N$. For $\Omega\subset\Lambda$ we have a canonical inclusion $\CZ^\Omega\subset\prod_{x\in\Omega} S$ and an induced inclusion $\CZ^\Omega_Q\subset \prod_{x\in\Omega} Q$.

\begin{lemma}[{\cite[Lemma 3.1]{Fie05}}] For a finite subset $\Omega$ of $\Lambda$ the inclusion above is a bijection, i.e.~
$\CZ^\Omega_Q=\bigoplus_{x\in\Omega} Q$.
\end{lemma} 

Choose $M\in\CZ\catmod^f$ and let $\Omega\subset\Lambda$ be finite such that the action of $\CZ$ on $M$ factors over $\CZ^\Omega$. Then $\CZ^\Omega_Q=\bigoplus_{x\in\Omega} Q$ acts on $ M_Q$ and the idempotents give a canonical decomposition
$$
 M_Q=\bigoplus_{x\in\Omega}  M_{Q}^x.
$$ 
We set $M_{Q}^x=0$ if $x\in\Lambda\setminus\Omega$. Since $M$ is supposed to be torsion free over $S$, there is a canonical injection $M\inj M_Q$ that  allows us to view each $m\in M$ as a $\Lambda$-tuple $(m_x)$ with $m_x\in M_{Q}^x$.

\begin{definition}
For $M\in\CZ\catmod^f$ and a subset $\Omega\subset\Lambda$ let
$M_{\Omega}$ be the biggest submodule of $M$ supported on $\Omega$, i.e.~ 
$$
M_\Omega:= M\cap\bigoplus_{x\in\Omega} M_{Q}^x,
$$
and let $M^{\Omega}$ be the biggest image of $M$ supported on $\Omega$, i.e.~ 
$$
M^\Omega:=\im\left(M\to M_Q\to\bigoplus_{x\in\Omega}M_{Q}^x\right),
$$ 
where $M_Q\to\bigoplus_{x\in\Omega}M_{Q}^x$ is the projection with kernel $\bigoplus_{x\in\Lambda\setminus \Omega}M_{Q}^x$.
\end{definition}
$M_\Omega$ and $M^\Omega$ are $\CZ$-modules, and for $\Omega\subset\Omega^\prime\subset\Lambda$ there is a canonical inclusion $M_\Omega\to M_{\Omega^\prime}$ and a canonical surjection $M^{\Omega^\prime}\to M^\Omega$. Write $M^{x,y,\dots,z}$ and $M_{x,y,\dots,z}$ for $M^{\{x,y,\dots,z\}}$ and $M_{\{x,y,\dots,z\}}$, respectively.

\subsection{A cofiltration associated to a partial order on $\Lambda$}
Suppose that $\leq$ is a partial order on the set $\Lambda$ such that $x,y\in\Lambda$ are comparable if $x=ty$ for some $t\in\CT$. If, for example, there is $x\in\Lambda$ such that $\CW^\prime=\Stab_\CW\, x$  is generated by a set of {\em simple} reflections $\CS^\prime\subset\CS$, then the Bruhat order on $\CW$ descends onto $\CW/\CW^\prime$ (compare the representatives of smallest length in each coset), and the map $w\mapsto w.x$ induces a bijection $\CW/\CW^\prime=\Lambda$. The induced partial order is of the above kind.  

Such a partial order on $\Lambda$ induces a cofiltration on $M$ as follows.  For $y\in\Lambda$ write $M^{\geq y}=M^{\{\geq y\}}$, where $\{\geq y\}=\{z\in\Lambda\mid z\geq y\}$. There is the natural quotient map $M^{\geq y}\to M^{\geq z}$ if $y\leq z$. Define $M^{>y}$ similarly and set
$$
M^{[y]}:=\ker \left(M^{\geq y}\to M^{>y}\right).
$$
The cofiltration and its subquotients are functorial, i.e.~ respected by every morphism in $\CZ\catmod^f$.

Let $M$ be a graded $S$-module. For $l\in\DZ$ denote by $M\{l\}$ the graded $S$-module that is obtained from  $M$ by shifting the grading by $l$, i.e.~ such that $M\{l\}_{\{n\}}=M_{\{n+l\}}$. A graded $S$-module $M$ is called {\em graded free} if there is an isomorphism $M\cong \bigoplus S\{k_i\}$ for finitely many $k_i\in\DZ$.

A subset $\Omega\subset\Lambda$ is said to be {\em upwardly closed} if $\Omega=\bigcup_{y\in\Omega}\{\geq y\}$. We can define a topology on $\Lambda$ with the upwardly closed sets as the open sets. 

\begin{definition}
We say that $M\in\CZ\catmod^f$ {\em admits a Verma flag} if, for any upwardly closed set $\Omega\subset\Lambda$, the module $M^\Omega$ is a graded free $S$-module. We denote by
$\CV=\CV(\Lambda,\leq)\subset \CZ\catmod^f$  the full subcategory of objects that admit a Verma flag.
\end{definition}

If $M$ admits a Verma flag, then $M^{[y]}$ is graded free for any $y\in \Lambda$. 

\subsection{An exact structure}\label{exact structure}
$\CV$ is an additive category, but it is not abelian in all but the trivial cases. As a substitute we will now provide $\CV$ with an exact structure in the sense of Quillen (cf.~ \cite{Qu}). An exact structure gives a meaning to notions such as exact functor or projective object. 

\begin{definition}
Let $A\to B\to C$ be a sequence in $\CV$. We say that it is {\em short exact} if for any upwardly closed set $\Omega\subset\Lambda$ the induced sequence
$$
0\to A^\Omega\to B^\Omega\to C^\Omega\to 0
$$
is a short exact sequence of $S$-modules (cf.~ Section \ref{secloc} for an example).
\end{definition}

\begin{proposition}[\cite{Fie05}, Theorem 4.1] The above definition gives an exact structure on  $\CV$.
\end{proposition}

In the following the term ``exact structure'' will refer to the exact structurejust defined. For the standard exact structure we reserve the term ``exact as a sequence of abelian groups''. Note that the functors $({\cdot})^{[x]}, ({\cdot})^{\geq x},\dots$ are exact functors from $\CV$ to $\CV$, i.e.~ they preserve short exact sequences.

Let $\Omega$ be an upwardly closed subset and choose a minimal element $y\in\Omega$. Let $K$ be the kernel of the map $M^\Omega\to M^{\Omega\setminus\{y\}}$. We have the following commutative diagram with rows that are exact as sequences of abelian groups (and even in our new sense):

\centerline{
\xymatrix{
0 \ar[r]& K\ar[r] \ar[d]& M^\Omega\ar[r] \ar[d]& M^{\Omega\setminus \{y\}}\ar[r]\ar[d] & 0 \\
0 \ar[r] & M^{[y]} \ar[r] &  M^{\geq y} \ar[r] &  M^{>y}\ar[r] &  0.  \\  
}
}
\noindent

The proof of the next lemma uses a realization of $M$ as a sheaf on a moment graph (cf.~ Section \ref{MomGra}).
\begin{lemma}[\cite{Fie05}, Lemma 4.3] For $M\in\CV$ the map $K\to M^{[y]}$ is an isomorphism.
\end{lemma}

The proposition means that we can linearize the cofiltration in the following way. Suppose that $\{\dots,x_{-1},x_0,x_1,\dots\}=\Lambda$ is an enumeration of $\Lambda$ such that $i<j$ if $x_i<x_j$. Let $\Omega_i=\bigcup_{j\geq i}\{\geq x_j\}$. Then each  $\Omega_i$ is upwardly closed and $\Omega_{i}=\Omega_{i+1}\cup\{x_i\}$ for all $i$. Let $M^{\geq i}:=M^{\Omega_i}$. This defines a cofiltration on $M$, indexed by the integers, such that $M^{[x_{i}]}=\ker(M^{\geq i}\to M^{\geq i+1})$. 

Let $x,y\in\Lambda$ and let $[x,y]=\{z\in\Lambda\mid x\leq z\leq y\}$ be the interval between $x$ and $y$. Note that $[x,y]=\{\geq x\}\setminus \bigcup_{z\geq x,y\not\geq z}\{\geq z\}$ is the difference of two upwardly closed sets. Hence we can view 
$$
M^{[x,y]}:=\ker\left(M^{\geq x}\to \sum_{z\geq x,y\not\geq z}M^{\geq z}\right)
$$
as a subquotient of the cofiltration.
The next lemma is easily proved by linearizing the cofiltration.

\begin{lemma}\label{def exactness} Let $A\to B\to C$ be a sequence in $\CV$. Then the following are equivalent:
\begin{enumerate}
\item $0\to A^{[x]}\to B^{[x]}\to C^{[x]}\to 0$ is an exact sequence of $S$-modules for any $x\in\Lambda$.
\item $0\to A^{[x,y]}\to B^{[x,y]}\to C^{[x,y]}\to 0$ is an exact sequence of $S$-modules for any interval $[x,y]\subset\Lambda$.
\item$0\to A^{\geq x}\to B^{\geq x}\to C^{\geq x}\to 0$ is an exact sequence of $S$-modules for any $x\in\Lambda$.
\item$0\to A^{\Omega}\to B^{\Omega}\to C^{\Omega}\to 0$ is an exact sequence of $S$-modules for any upwardly closed $\Omega\subset\Lambda$.
\end{enumerate}
\end{lemma}

\subsection{A duality}
For $M\in\CZ\catmod^f$ set $D(M)=\bigoplus_{i\in\DZ} \Hom_S^i(M,S)$, where $\Hom^i_S(M,S)=\Hom_S(M,S\{i\})$. We have $D\circ\{n\}=\{-n\}\circ D$, and $\CZ$ naturally acts on $D(M)$ by $z.f(m)=f(z.m)$ for $z\in\CZ, f\in D(M)$ and $m\in M$. Then $D(M)\in\CZ\catmod^f$, so $D$ is a duality (i.e.~ $D^2=\id$) on the subcategory of objects that are reflexive over $S$. Moreover, $D$ inverts the canonical filtration. If $M\in\CV$, then each quotient, hence also each kernel is graded free. So $D$ induces an equivalence of exact categories
$$
 \CV(\Lambda)^{opp}\cong\CV(\Lambda^\circ),
$$
where $\Lambda^\circ$ is the set $\Lambda$ with the partial order reversed, and where $\CA^{opp}$ denotes the opposite category of a category $\CA$.

\section{Sheaves on moment graphs}\label{MomGra}
In this section we introduce the local definition of the objects of $\CV$ as sheaves on a moment graph.
 \subsection{Moment graphs}
Let $k$ be a field and $W$ a vector space over $k$. 
A {\em $W$-moment graph} $\CG=(\CV,\CE,\leq,l)$ is given by (cf.~ \cite{BMP}) 
\begin{itemize}
\item a (not necessarily finite) graph $(\CV,\CE)$ with a set of vertices $\CV$ and a set of edges $\CE$,
\item a partial order $\leq$ on $\CV$ such that $x$ and $y$ are comparable if they are linked by an edge, and
\item a labelling $l\colon\CE\to \DP\, W$ by one-dimensional subspaces of $W$.
\end{itemize}
We can think of the order as giving each edge a direction. We write $E\colon x\to y$ for an edge $E$ with endpoints $x$ and $y$ such that $x\leq y$, and $E\colon x\xrightarrow{\alpha}y$, where $\alpha\in W$ is a generator of $l(E)$ if we want to denote the label. We write $ x\linie y$ or $x\stackrel{\alpha}\linie y$ if we want to ignore the order.

\subsection{Sheaves on moment graphs}
Let $S=S_k(W)$ be the symmetric algebra of $W$ with a grading given by $\deg W=2$. A {\em sheaf ${\SM}=(\{\SM^x\}, \{\SM^E\}, \{\rho_{x,E}\})$ on $\CG$} is given by 
\begin{itemize}
\item a graded $S$-module $\SM^x$ for each vertex $x\in\CV$,
\item a graded $S$-module $\SM^E$  for each edge $E\in\CE$ such that $l(E)\cdot \SM^E=0$,  
\item a map $\rho_{x,E}\colon \SM^x\to \SM^E$ of graded $S$-modules if $x$ is a vertex of $E$.  \end{itemize}
In addition we assume that each $\SM^x$ is torsion free and finitely generated over $S$, and that $\SM^x$ is non-zero for only finitely many $x$. A morphism $\SM\to\SN$ of sheaves on $\CG$ is given by maps $\SM^x\to\SN^x$ and $\SM^E\to\SN^E$ that are compatible with the maps $\rho$. Let $\Sh(\CG)$ be the resulting category.

Let $\SM$ be a sheaf on $\CG$. 
The {\em space of global sections of $\SM$} is
$$
\Gamma(\SM):=
\left\{ 
(m_x)\in\prod_{x\in \CV}  \SM^x\left| \,
\begin{matrix}
\rho_{x,E}(m_x)=\rho_{y,E}(m_{y}) \\
\text{for any edge $E\colon x\linie y$}
\end{matrix}
\right.
\right\}.
$$
For a full subgraph $\Omega\subset\CG$ we write $\SM|_\Omega$ for the restriction of $\SM$ to $\Omega$ and $\Gamma(\Omega,\SM)=\Gamma(\SM|_\Omega)$ for the sections of $\SM$ over $\Omega$. 

The {\em structure algebra} of $\CG$ is defined by 
$$
\CZ=\CZ(\CG)=\left\{ 
(z_x)\in\prod_{x\in \CV}  S\left| \,
\begin{matrix}
z_x\equiv z_y\mod \alpha  \\
\text{for any edge $E\colon x\stackrel{\alpha}\llinie y$}
\end{matrix}
\right.
\right\}.
$$
The global sections of any sheaf $\SM$ on $\CG$ naturally form a $\CZ$-module by pointwise multiplication.

Let $x$ be a vertex of $\CG$ and let $\CE^{\delta x}=\{E\colon x\to y\}$ be the set of edges {\em starting} at $x$. Define $\SM^{[x]}\subset\SM^x$ as the kernel of the map $\bigoplus_{E\in\CE^{\delta x}} \rho_{x,E}\colon \SM^x\to\bigoplus \SM^E$.

\subsection{The moment graph associated to a reflection representation and localization}
Let $V$ be a reflection faithful representation of $(\CW,\CS)$ and let $\Lambda\subset V$ be an orbit endowed with a partial order as before. The associated moment graph $\CG=\CG(\Lambda)$ is defined as follows. It is a moment graph over $V^\ast$. Its partially ordered set of vertices is $\Lambda$, and two vertices $x,y\in\Lambda$, $x\neq y$, are connected by an edge if there is a reflection $t\in\CT$ such that $y=tx$. The edge is labelled by $k\cdot \alpha_t\in\DP\, V^\ast$. Note that by Lemma \ref{rfrep}, (2), there are no double edges, i.e.~ two vertices are connected by at most one edge. 

Note that $\CZ(\CG)$ is just the structure algebra $\CZ(\Lambda)$ that was defined earlier. 
In \cite{Fie05} we showed that one can view each object of $\CV$ as a sheaf on $\CG$ via a {\em localization functor} $\Loc\colon \CZ\catmod^f\to\Sh(\CG)$ that is left adjoint to $\Gamma$. It induces a full embedding $\CV\stackrel{\sim}\to\Loc(\CV)\subset\Sh(\CG)$ with left inverse $\Gamma$. 

The definition of $\Loc$ is as follows. Let $M\in\CZ\catmod^f$. For a vertex $x$ we have $\Loc(M)^x:= M^x$ and for an edge $E\colon x\linie y$ the module $\Loc(M)^E$, together with the maps $\rho_{x,E}$ and $\rho_{y,E}$, is defined as the bottom right corner in the following push out diagram:

\centerline{
\xymatrix{
M^{x,y}\ar[d]\ar[r] & M^x\ar[d]^{\rho_{x,E}} \\
M^y\ar[r]^{\rho_{y,E}}&\Loc(M)^E.
}
}
\noindent
(Note that in our situation and in the notation of \cite{Fie05} we have $M(E)=M^{x,y}$ since the map $\CZ\to\CZ(E)$ is surjective by Lemma \ref{loccen}). 
\begin{lemma}[\cite{Fie05}, Lemma 4.8]\label{loccar} Let $M\in\CZ\catmod^f$. Then $M\in\CV$ if and only if for $\SM=\Loc(M)$ the adjunction morphism $M\to\Gamma(\SM)$ is an isomorphism and the following holds:
\begin{enumerate}
\item $\SM$ is flabby, i.e.~ for each upwardly closed subset $\Omega\subset\Lambda$ the canonical map 
$$
\Gamma(\SM)\to\Gamma(\SM|_{\Omega})
$$
is surjective.
\item For each $x\in\Lambda$ the $S$-module $\SM^{[x]}$ is graded free. 
\end{enumerate}
\end{lemma}

\subsection{The local structure of modules with a Verma flag}\label{secloc}
For $x\in\Lambda$ let $M(x)$ be the $\CZ$-module that is free of rank one over $S$ and on which $(z_v)\in\CZ$ acts by multiplication with $z_x$. Then  $M(x)$ is supported on $\{x\}$ and admits a Verma flag. Note that each $M\in\CV$ is an extension of objects of the form $M(x)\{k\}$ with $x\in\Lambda$ and $k\in\DZ$.

Let $x,y\in\Lambda$ such that $y=tx\neq x$ for some $t\in\CT$. Denote by $E\colon x\stackrel{\alpha_t}\linie y$ the corresponding edge and let
$$
\CZ(E)=\left\{(z_x,z_y)\in S\oplus S\mid z_x\equiv z_y\mod\alpha_t\right\}
$$
be the {\em local structure algebra}. Projection onto the coordinates $x$ and $y$ gives an algebra map $\CZ\to\CZ(E)$, so we can view $\CZ(E)$ as a $\CZ$-module. It admits a Verma flag and we denote it  by $P(x,y)$. Suppose that $x<y$. Then $P(x,y)_x=P(x,y)^{[x]}=S\cdot (\alpha_t,0)=S\{-2\}$ and $P(x,y)^{[y]}=P(x,y)^y=S$. As an example of the exact structure defined earlier, note that the sequence
$$
P(x,y)_x\to P(x,y)\to P(x,y)^y
$$
is exact, while
$$
P(x,y)_y\to P(x,y)\to P(x,y)^x
$$
is not.

\begin{remark}\label{locstruc} By Lemma \ref{loccen} the map  $\CZ\to\CZ(E)$ is surjective. Thus a $\CZ(E)$-module is the same as a $\CZ$-module supported on $\{x,y\}$. 
\end{remark}

\begin{lemma}\label{local structure} Let $M$ be a $\CZ(E)$-module that is graded free over $S$ of finite rank. Then $M$ is isomorphic to a direct sum of shifted copies of $M(x)$, $M(y)$ and $P(x,y)$. 
\end{lemma}
\begin{proof} For a graded $S$-module $N$ and $k\in\DZ$ let $N_{\{\leq k\}}=\sum_{i\leq k} S\cdot N_{\{i\}}$ be the submodule generated in degrees $\leq k$. 

We can assume that $M$ is evenly graded. Suppose that we have already shown that the submodule $X=\CZ(E)\cdot M_{\{\leq 2(n-1)\}}$ is isomorphic to a direct sum of the alleged kind, and let $P_{n-1}\subset X$ be the direct sum of the summands isomorphic to $P(x,y)\{-2(n-1)\}$. Then $X=M_{\{\leq 2(n-1)\}}\oplus (0,\alpha_t)\cdot P_{n-1}$, where $(0,\alpha_t)\cdot P_{n-1}$ is isomorphic to a direct sum of copies of $M(y)\{-2n\}$. 

Choose a $k$-complement $A$ of $X\cap (M_x)_{\{2n\}}$ in $(M_x)_{\{2n\}}$, a $k$-complement $B$ of $X\cap (M_y)_{\{2n\}}$ in $(M_y)_{\{2n\}}$ and a $k$-complement $C$ of $X\cap M_{\{2n\}}+(M_x)_{\{2n\}}+(M_y)_{\{2n\}}$ in $M_{\{2n\}}$.  By construction,
$$
\CZ(E)\cdot M_{\{\leq 2n\}}=X\oplus\CZ(E)\cdot A\oplus\CZ(E)\cdot B\oplus\CZ(E)\cdot C
$$
and $\CZ(E)\cdot A$ is isomorphic to a direct sum of copies of $M(x)\{-2n\}$,  $\CZ(E)\cdot B$ is isomorphic to a direct sum of copies of $M(y)\{-2n\}$ and $\CZ(E)\cdot C$ is isomorphic to a direct sum of copies of $P(x,y)\{-2n\}$.
\end{proof}

\subsection{The connection to representation theory} We want to explain the representation theoretic content of $\CV$. The following will not be used in the sequel. More details can be found in \cite{Fie05}.

Let $\fg$ be a symmetrizable, complex Kac-Moody algebra and $\fh\subset\fb\subset\fg$ its Cartan and Borel subalgebras. Let $S=S(\fh)$ and $R=S_{(0)}$ be the localization at the point $0\in\fhd$. Let $\tau\colon \fh\to R$ be the canonical map. For each $\fh$-$R$-bimodule $M$ and each $\lambda\in\fhd$ let $M_\lambda=\{m\in M\mid H.m=(\lambda(H)+\tau(H))m \text{ for all $H\in\fh$}\}$ be its {\em weight space of weight $\lambda$}. Define $\CO_R$ as the subcategory of the category of $\fg$-$R$-bimodules that consists of all modules $M$ such that $M=\bigoplus_{\lambda\in\fhd}M_\lambda$ and such that for each $m\in M$, the $\fb$-$R$-submodule generated by $m$ is finitely generated over $R$. We call it the {\em equivariant} category $\CO$, since it is related to the torus-equivariant topology on the associated flag variety (cf.~ \cite{Fie05}).

For any $\lambda\in\fhd$ let $R_\lambda$ be the $\fh$-$R$-bimodule that is free of rank one over $R$ and on which $\fh$ acts via the character $\lambda\cdot 1+\tau\colon \fh\to R$. By letting $[\fb,\fb]$ act trivially we extend the action to the Lie algebra $\fb$, and induction yields
$$
M_R(\lambda):=U(\fg)\otimes_{U(\fb)}R_\lambda,
$$ 
the {\em equivariant Verma-module with highest weight $\lambda$}.
It is an object of $\CO_R$. We define $\CM_R\subset\CO_R$ as the subcategory of all objects that admit a finite filtration with subquotients isomorphic to Verma modules. It is called the category of equivariant modules that admit a Verma flag. 

Choose an indecomposable block $\CO_{R,\Lambda}$ of $\CO_R$. Then we can view $\Lambda$ as the set of highest weights of the Verma modules in $\CO_{R,\Lambda}$, hence as a subset of the dual space $\fhd$ of $\fh$. If $\Lambda$ lies outside the critical hyperplanes, it is an orbit of its integral Coxeter system $(\CW_\Lambda,\CS_\Lambda)$ and there is a natural order on $\Lambda$ that is given by $x\leq y$ if and only if $y-x$ is a sum of positive roots. By \cite[Proposition 2.1]{Soe04}, $\fh$ and $\fhd$ are reflection faithful representations.

Let $\CG(\Lambda)$ be the moment graph associated to the orbit $\Lambda$. 
In all the constructions and definitions so far we can replace the symmetric algebra $S$ by its localization $R$, and we denote the resulting objects with an index $R$. To $\CG=\CG(\Lambda)$ we associate the $R$-linear exact category $\CV_R(\Lambda)$. 

Let $\CM_{R,\Lambda}\subset\CO_{R,\Lambda}$ be the corresponding block of $\CM$. It inherits an exact structure from the abelian category $\CO_R$.

\begin{theorem}[\cite{Fie05}, Theorem 7.1] There is an equivalence $\CM_{R,\Lambda}\cong\CV_R(\Lambda)$  of exact categories. 
\end{theorem}

Since all projective objects of $\CO_{R,\Lambda}$ lie in $\CM_{R,\Lambda}$, $\CV_R(\Lambda)$ captures all the categorical information of $\CO_{R,\Lambda}$. 
\section{$S$-bimodules}\label{S-bimod}
 Soergel associated in \cite{Soe04} a category $\CF_\Nabla$ of $S$-bimodules to any reflection faithful representation $V$.  We want to review his construction and show that there is an equivalence $\CF_\Nabla\cong \CV(\Lambda)$ for each regular orbit $\Lambda\subset V$. Moreover, Soergel's subcategory $\SB$ of {\em special bimodules} turns out to correspond to the subcategory of projective objects with respect to the exact structure defined in \ref{exact structure} (cf.~ Theorem \ref{BequivP}).

Let $V$ be a reflection faithful representation of $(\CW,\CS)$ and let $S=S(V^\ast)$. Suppose that $k$ is infinite and identify $S$ with the algebra of regular functions on $V$.
Let $S\bimod S$ be the category of $S$-bimodules. We view such a bimodule  as a quasi-coherent sheaf on the variety $V\newtimes  V$. For $x\in\CW$ let $\Gr(x)=\{(x^{-1}.v,v)\in V\newtimes  V\mid v\in V\}$ be the twisted diagonal (in \cite{Soe04} this variety is denoted by $\Gr(x^{-1})$, we have chosen the inverse to simplify notations).  An important property of a reflection faithful representation is the following.

\begin{lemma}[{\cite[Bemerkung 3.2]{Soe04}}]\label{int of hyperplanes} Let $x,y\in\CW$. Then $\Gr(x)\cap\Gr(y)$ is of codimension one in $\Gr(x)$ if and only if $y=tx$ for a reflection $t\in\CT$. 
\end{lemma}

Let $S(x)$ be the space of regular functions on $\Gr(x)$. It is naturally an $S$-bimodule. The projection $pr_2\colon \Gr(x)\subset V\newtimes  V\to V$ onto the second component is an isomorphism of varieties. This identifies $S(x)$ with $S$ in such a way that the element $f\otimes g\in S\otimes_kS$ acts on $S(x)$ by multiplication with $f^x\cdot g\in S$, where $f^x\in S$ is the function that maps $v\in V$ to $f(x^{-1}.v)$.

For a subset $A$ of $\CW$ set $\Gr(A)=\bigcup_{x\in A}\Gr(x)$ and for $M\in S\bimod S$ let $M_A\subset M$ be the submodule of sections supported on $\Gr(A)$. In particular, for $i\geq0$,  let $M_{\leq i}\subset M$ be the submodule supported on $\Gr(\{x\mid l(x)\leq i\})$, where $l\colon \CW\to\DN$ is the length function associated to $\CS$.

\begin{definition}[\cite{Soe04}] Let $\CF_\Nabla\subset S\bimod S$ be the full subcategory of objects $M$ that are supported on $\Gr(A)$ for some {\em finite} subset $A$ of $\CW$ and such that for each $i\geq 0$, $M_{\leq i+1}/M_{\leq i}$ is isomorphic to a finite direct sum of graded bimodules of the form $S(x)\{k\}$ with $l(x)=i$ and $k\in\DZ$.
\end{definition}

Choose a {\em regular} orbit $\Lambda\subset V$, i.e.~ an orbit whose stabilizers are trivial. The choice of an element $v\in\Lambda$ gives an identification $\CW\cong\Lambda$ that we fix. The Bruhat order induces a partial order on $\Lambda$. Let $\CZ=\CZ(\Lambda)$ and $\CV=\CV(\Lambda)$.

\begin{theorem}\label{eqcat} There is an equivalence $\CF_\Nabla\cong \CV$ of categories.
\end{theorem}

We will see in the course of the proof of the theorem that $\CZ$ can be canonically identified with the regular functions on the union $\cup_{x\in\CW}\Gr(x)$, at least for finite $\CW$.

Before we prove the theorem we need some preparations.
For $x\in\CW$ and $t\in\CT$ consider the regular functions $S(x,tx)$ on the union $\Gr(x)\cup\Gr(tx)$. Let $S(x,tx)\to S(x)$ and $S(x,tx)\to S(tx)$ be the restrictions to $\Gr(x)$ and $\Gr(tx)$. Then $S(x,tx)\to S(x)\oplus S(tx)$ is injective and identifies $S(x,tx)$ with the set of pairs $(f_x, f_{tx})\in S(x)\oplus S(tx)$  that agree on the intersection $\Gr(x)\cap\Gr(tx)$. We want to  generalize this representation.

Choose $M\in\CF_\Nabla$ and denote by $M^x$ the restriction  of $M$ to $\Gr(x)$. If $M$ is supported on $\Gr(A)$ for some subset $A$ of $\CW$, then the canonical map $M\to \bigoplus_{x\in A} M^x$ is injective. Similarly, let $M^{x\cap y}$ be the restriction to $\Gr(x)\cap\Gr(y)$. Consider the restriction maps $\rho_{x,x\cap y}\colon M^x\to M^{x\cap y}$ and $\rho_{y,x\cap y}\colon M^y\to M^{x\cap y}$. We consider these objects as $S$-modules using the right action of $S$. The following proposition is analogous to the localization theorem in equivariant topology.
\begin{proposition}\label{loc bimod} The inclusion $M\inj \bigoplus_{x\in A} M^x$ gives an identification
$$
M=
\left\{
(m_x)\in\bigoplus_{x\in A} M^x 
\left|
\,\begin{matrix} \rho_{x,x\cap tx}(m_x)=\rho_{tx,x\cap tx}(m_{tx}) \\
\text{for all $t\in\CT$, $x,tx\in A$}
\end{matrix}
\right.
\right\}.
$$ 
\end{proposition}
\begin{proof}
For a graded prime ideal $\fp\subset S$ let $S_\fp$ be the graded localization at $\fp$ and for $N\in S\bimod S$ define $N_\fp:=N\otimes_S S_\fp\in S\bimod S_\fp$.  Recall that we denote by $\alpha_t\in V^\ast$  a generator of $(V^\ast)^{-t}$, i.e.~ an equation of the reflection hyperplane $V^t$.

\begin{lemma}\label{ext of loc Verma} Let $\fp\subset S$ be a graded prime ideal of height one.
\begin{enumerate}
\item If $\alpha_t\in\fp$ for some $t\in\CT$, then $\fp=S\cdot \alpha_t$ and $M_\fp$ is isomorphic to a direct sum of shifted copies of $S(x)_\fp$  and $S({x,tx})_\fp$ for $x\in\CW$.
\item If $\alpha_t\not\in \fp$ for all $t\in\CT$, then  $M_\fp$ is isomorphic to a direct sum of shifted copies of $S(x)_\fp$ for $x\in\CW$.
\end{enumerate}
\end{lemma}

\begin{proof}
$M$ is an extension of modules of the form $S(x)$, $x\in\CW$. By \cite{Soe04}, Lemma 5.8, $S(x)$ and $S(y)$ do not extend unless $\Gr(x)\cap\Gr(y)$ is of codimension one in $\Gr(x)$. By Lemma \ref{int of hyperplanes} this is the case if and only if  $y=tx$ for some $t\in\CT$. Moreover, again by Lemma 5.8 in \cite{Soe04}, each extension splits after inverting $\alpha_t$, and the class of extensions of $S(x)_{S\cdot\alpha_t}$ and $S(tx)_{S\cdot \alpha_t}$ is generated by $S(x,tx)_{S\cdot \alpha_t}$.
\end{proof}

Now we can prove the proposition. Denote by $\widehat M$ the space on the right hand side of the alleged identity. Since the right action of $S$ on $M$ is graded free, $M$ is the intersection of all its localizations $M_\fp$ for graded prime ideals $\fp\subset S$ of height one. So we only have to show that for each such $\fp$ the localization $M_\fp\to\widehat M_\fp$ is an isomorphism. Lemma \ref{ext of loc Verma} reduces this to the cases $M_\fp=S(x)_\fp$ and $M_\fp=S(x,tx)_\fp$. The first case is clear and the second is the example we started with. 
\end{proof}

\begin{proof}[Proof of Theorem \ref{eqcat}] We will define a functor $G\colon \CF_\Nabla\to\CV$ by representing a module $M\in\CF_\Nabla$ as the global sections of a sheaf on the moment graph associated to $\Lambda$, and a functor $F\colon \CV\to\CF_\Nabla$ as a restriction along a map $S\otimes_k S\to \CZ$.

Let $M\in\CF_\Nabla$. Define a sheaf $\SM$ on $\CG$ by setting $\SM^x:=M^x$, $\SM^{E}:=M^{x\cap tx}$ for any edge $E\colon x\linie tx$, and $\rho_{x,E}=\rho_{x,x\cap tx}\colon \SM^x\to\SM^E$. Note that the projection $pr_2\colon V\newtimes V\to V$ onto the second factor induces an isomorphism $\Gr(x)\cap\Gr(y)\cong V^{yx^{-1}}$. For $y=tx$ we deduce that $M^{x\cap tx}$ is annihilated by the right action of $\alpha_t$. Hence we have defined a sheaf $\SM$ on $\CG$. By Proposition \ref{loc bimod} we have $M=\Gamma(\SM)$, hence there is a canonical action of $\CZ$ on $M$. This gives a functor $G\colon \CF_\Nabla\to \CZ\catmod$. 

Let $\Omega\subset\CW$ be upwardly closed and let $M^\Omega$ be the restriction to $\Gr(\Omega)$. It follows from Lemma 6.3 in \cite{Soe04} that $M^\Omega$ lies in $\CF_\Nabla$ as well. Clearly $G(M)^\Omega=G(M^\Omega)$. In particular, $G(M)^\Omega$ is a graded free $S$-module, hence $G(M)\in\CV$ and we get a functor $G\colon\CF_\Nabla\to\CV$.

We will now construct the functor $F$. Note that $\CW$ acts on $S$ by $w.f(v)=f(w^{-1}.v)$ for $w\in\CW$, $f\in S$ and $v\in V$. 
Recall the homomorphism $\sigma\colon V^\ast\to \CZ$, $\sigma(\lambda)_w=w.\lambda$,  defined in Remark \ref{sigmareg}. 
Let $\tau\colon S\to\CZ$, $f\mapsto f\cdot 1$, be the $S$-structure on $\CZ$. We get a homomorphism 
$$
\sigma\otimes \tau\colon S\otimes_k S\to \CZ
$$
of algebras and a corresponding restriction functor $F\colon\CV\to S\bimod S$.

Let $x\in\CW$. Recall that we identified the regular functions on $\Gr(x)$ with the regular functions on $V$ by means of the isomorphism $pr_2\colon \Gr(x)\stackrel{\sim}\to V$. 
Let $i_x\colon S\otimes_k S\to S$, $i_x(f\otimes g)=f^x\cdot g$, be the map corresponding to the inclusion $\Gr(x)\to V\newtimes  V$. Then following diagram commutes:

\centerline{
\xymatrix{
S\otimes_k S \ar[dr]_{i_x}\ar[rr]^{\sigma\otimes \tau} & & \CZ \ar[dl]^{\eval_x} \\
& S &  
}
}
\noindent
where $\eval_x((z_v)_v)=z_x$. Hence $F$ respects the canonical filtrations, i.e.~ $F(M)_{\leq i}= F(M_{\{l(\cdot)\leq i\}})$ where $\{l(\cdot)\leq i\}=\{x\in\CW\mid l(x)\leq i\}$. From linearizing the cofiltration of $M$ we deduce that $M_{\{l(\cdot)\leq i\}}$ is an extension of modules of the form $M(x)\{k\}$ with $l(x)\leq i$ and $k\in\DZ$. So $M_{\{l(\cdot)\leq i\}}/M_{\{l(\cdot)\leq i-1\}}$ is an extension of modules of the form $M(x)\{k\}$ with $l(x)= i$. By Lemma \ref{loccar}, $M$ is the space of  global sections of a sheaf on $\CG$. But two vertices of the same length are not connected, so the extension splits. From $F(M(x)\{k\})\cong S(x)\{k\}$ we deduce $F(M)\in\CF_\Nabla$. 

It is easily verified that  $F$ and $G$ are mutually inverse equivalences $\CF_\Nabla\cong\CV$ of categories.
\end{proof}

\section{Translation functors}
Suppose that $V$ is a reflection faithful representation of $\CW$ and that $\Lambda=\CW.v\subset V$ is a regular orbit, i.e.~ $\Stab_\CW(v)=\{e\}$. Choose a simple reflection  $s\in\CS$. An {\em $s$-subregular} orbit is an orbit $\Lambda^\prime=\CW.v^\prime$ of an element $v^\prime$ such that $\Stab_\CW(v^\prime)=\{e,s\}$. Fix such a $\Lambda^\prime$. We identify $\Lambda$ with $\CW$ and $\Lambda^\prime$ with $\CW/s$. Let $\CZ=\CZ(\Lambda)$ and $\CZ^\prime=\CZ(\Lambda^\prime)$ be the associated algebras, and $\CV=\CV(\Lambda)$ and $\CV^\prime=\CV(\Lambda^\prime)$ the associated exact categories. For $x\in\CW$ write $\overline{x}\in\CW/s$ for its canonical image. More generally, we write  $\overline{\CI}\subset\CW/s$ for the image of $\CI\subset\CW$.

\subsection{An action of $\CW$ on $\CZ$}\label{WonZ}

Let $\CW$ act on the left of $\prod_{w\in\CW} S$ by switching coordinates, i.e.~  $x.(z_w)=(z^\prime_w)$ with $z^\prime_w=z_{wx}$. Then
$\CZ\subset \prod_{w\in\CW} S$ is $\CW$-stable. The diagonal inclusion $\prod_{\overline{w}\in\CW/s}S\to\prod_{w\in\CW}S$, $(z_{\overline{w}})\mapsto (z^\prime_w)$ with $z^\prime_{ws}=z^\prime_w=z_{\overline{w}}$, induces a map $\CZ^\prime\to\CZ$ that identifies $\CZ^\prime$ with the $s$-invariants $\CZ^s\subset\CZ$.  This allows us to define 
the functor $\out\colon\CZ^\prime\catmod\to\CZ\catmod$ of {\em translation out of the wall} by 
$$
\out(N):=\CZ\otimes_{\CZ^\prime}N,
$$ 
and the {\em translation on the wall}, $\on\colon\CZ\catmod\to\CZ^\prime\catmod$, as the restriction, the right adjoint to $\out$. 
In the remainder of this section we will study the properties of $\on$ and $\out$.

Recall the map $\sigma\colon V^\ast\to\CZ$, defined in Remark \ref{sigmareg}. Let $c^s:=\sigma(\alpha_s)$. Then $c^s$ is $s$-anti-invariant, i.e.~ $s.c^s=-c^s$, since $s(\alpha_s)=-\alpha_s$ (recall that $\ch\, k\neq 2$). 

\begin{lemma}\label{dec of Z} There is a decomposition
$\CZ=\CZ^s\oplus c^s\cdot\CZ^s$
of $\CZ^s$-modules.
\end{lemma}

\begin{proof} Note that $\CZ=\CZ^s\oplus\CZ^{-s}$, where $\CZ^{-s}$ is the set of $s$-anti-invariant elements. Hence we have to show that each element in $\CZ^{-s}$ is divisible by $c^s$ in $\CZ$. So let $z\in\CZ^{-s}$. For each $w\in\CW$, $z_w=-z_{ws}$ and, by definition, $z_w\equiv z_{ws}\mod\alpha_{wsw^{-1}}$. Hence $z_w$ is divisible by $\alpha_{wsw^{-1}}$. Now $c^s_w=w(\alpha_s)$ is an equation of the reflection hyperplane of $wsw^{-1}$ as well, hence $c^s_w$ is a multiple of $\alpha_{wsw^{-1}}$. So $z^\prime_w=z_w/c^s_w\in S$. Set $z^\prime=(z^\prime_w)$. It is clearly $s$-invariant, so we only have to show that $z^\prime\in\CZ$.

Choose $w\in\CW$ and  $t\in\CT$. Then 
$(c^s_{tw}c^s_w)(z^\prime_{tw}-z^\prime_w)  =  c^s_{w}\cdot z_{tw}-c^s_{tw}\cdot z_w \equiv 0\mod \alpha_t$,
since $c^s_{tw}\equiv c^s_w$  and  $z_{tw}\equiv z_w \mod\alpha_t$. If $t\neq wsw^{-1}$, then $c^s_{tw}=tw(\alpha_s)\not\equiv 0\mod \alpha_t$ and $c^s_{w}=w(\alpha_s)\not\equiv 0\mod \alpha_t$, since $\alpha_t$ and $\alpha_{t^{\prime}}$ are linearly independent for $t\neq t^\prime$. Hence $z^\prime_{tw}-z^\prime_w\equiv 0\mod \alpha_t$. If $t=wsw^{-1}$, then $z^\prime_{tw}-z^\prime_w=0$. Hence $z^\prime_{tw}-z^\prime_w\equiv 0\mod \alpha_t$ for all $t\in\CT$ and $w\in\CW$, so  $z^\prime\in \CZ$. 
\end{proof}

\subsection{Adjointness}
We want to show that the translation functors $\on$ and $\out$ are biadjoint.
\begin{proposition}[cf.~ \cite{Soe04}, Proposition 5.10]\label{adj} The functors from $\CZ^s\catmod^f$ to $\CZ\catmod^f$ that are defined by $ M\mapsto \CZ\{2\}\otimes_{\CZ^s} M$ and $ M\mapsto\Hom_{\CZ^s}(\CZ, M)$, are isomorphic. In particular, $\on$ and $\out$ are biadjoint up to a shift, i.e.~ there are isomorphisms $\Hom(\out\cdot,\cdot)\cong\Hom(\cdot,\on\cdot)$ and $\Hom(\on\cdot,\cdot)\cong\Hom(\cdot,\{2\}\circ\out\cdot)$ of bifunctors.
\end{proposition}

\begin{proof}[Proof (due to Soergel)] Let $1^\ast,c^{s,\ast}\in\Hom_{\CZ^s}(\CZ,\CZ^s)$ be the $\CZ^s$-basis dual to $1,c^s$. Note that $\deg 1^\ast=\deg 1=0$ and $\deg c^s=2$, $\deg c^{s,\ast}=-2$. Hence $1\mapsto c^{s,\ast}$ and $c^s\mapsto 1^\ast$ gives a $\CZ^s$-isomorphism $\CZ\{2\}\cong\Hom_{\CZ^s}(\CZ,\CZ^s)$.  We have a functorial isomorphism $\Hom_{\CZ^s}(\CZ, M)=\Hom_{\CZ^s}(\CZ,\CZ^s)\otimes_{\CZ^s} M$ since $\CZ$ is free over $\CZ^s$. Hence the first claim.

Biadjointness follows from the fact that the functor $\CZ\otimes_{\CZ^s}\cdot$ is left adjoint, while  $\Hom_{\CZ^s}(\CZ,\cdot)$ is right adjoint to the restriction functor. 
\end{proof}

\subsection{The coordinates of $\on$ and $\out$}

\begin{proposition}\label{coor} Let $M\in\CZ\catmod^f$ and $N\in\CZ^\prime\catmod^f$. Then the following holds:
\begin{enumerate}
\item For $\CI\subset\CW$ such that $\CI\cdot s=\CI$ we have
$$
(\on M)^{\overline{\CI}}=\on(M^\CI),\quad (\out N)^{\CI}=\out(N^{\overline{\CI}})
$$
 and
$$
(\on M)_{\overline{\CI}}=\on(M_\CI),\quad (\out N)_{\CI}=\out(N_{\overline{\CI}}).
$$
\item  For $y\in\CW$ we have 
$$
(\on M)^{[\overline{y}]}  = \on(   M^{[y,ys]}),\quad (\out N)^{[y,ys]} = \out( N^{[\overline{y}]}).
$$ 
If $ys<y$, then 
$$
(\out N)^{[y]}\cong N^{[\overline{y}]}\cong (\out N)^{[ys]}\{2\}.
$$ 
\end{enumerate}
\end{proposition}
\begin{proof}
Suppose that the action of $\CZ$ on $M$ factors over $\CZ^\Omega$ with $\Omega\subset\CW$ finite. Then the action of $\CZ^\prime$ on  $\on M$ factors over $(\CZ^{\prime})^{\ol\Omega}$.  After enlarging $\Omega$ if necessary we can assume that $\Omega$ is right $s$-invariant, i.e.~ $\Omega\cdot s=\Omega$, and we can identify  $(\CZ^\prime)^{\ol \Omega}_Q$ with $(\bigoplus_{w\in\Omega} Q)^s$, the space of $s$-invariants in $\bigoplus_{w\in\Omega} Q$. So  $(\on M)_{Q}^{\overline{y}}= M_{Q}^y\oplus M_{Q}^{ys}$ for any $y\in\CW$. Hence $(\on M)_{\overline{\CI}}=\on(M_\CI)$ and $(\on M)^{\overline{\CI}}=\on(M^\CI)$ for right $s$-invariant $\CI\subset\CW$. 

Suppose that $ys<y$. Since the sequence
$$
0\to\on(M^{[y,ys]})\to \on(M^{\geq ys})\to \on(M^{> ys,\neq y})\to 0  
$$
is exact (as a sequence of abelian groups), and since $\on(M^{\geq ys})=(\on M)^{\geq \overline y}$ and $\on(M^{> ys,\neq y})= (\on M)^{>\overline{y}}$,  we have 
$$
(\on M)^{[\overline{y}]}=\on (M^{[y,ys]}).
$$ 

Using the above and Proposition \ref{adj}  we have, for $s$-invariant $\CI$,
\begin{eqnarray*}
\Hom((\out N)^\CI,M) & = & \Hom(\out N,M_{\CI}) \\
& = & \Hom(N,\on(M_{\CI}))\\
& = &  \Hom(N,(\on M)_{\ol\CI})\\
& = & \Hom(N^{\ol\CI},\on M)\\
& = & \Hom(\out(N^{\ol\CI}),M),
\end{eqnarray*}
hence $(\out N)^\CI=\out(N^{\ol\CI})$. Similarly we show that $(\out N)_\CI=\out(N_{\ol\CI})$. 

By Lemma \ref{dec of Z} the sequence
$$
0\to \out(N^{[\overline{y}]})\to\out(N^{\geq \overline{y}})\to\out(N^{>\overline{y}})\to 0
$$
is exact (as a sequence of abelian groups), hence $(\out N)^{[y,ys]}=\out(N^{[\overline{y}]})$. Moreover, $\out(N^{[\overline{y}]})=\CZ^{y,ys}\otimes_S N^{[\overline{y}]}$, where $\CZ^{y,ys}=\{(z_y,z_{ys})\in S\oplus S\mid z_y\equiv z_{ys}\mod\alpha_{ysy^{-1}}\}$.  It follows that $\out(N)^{[y]}\cong N^{[\overline{y}]}\cong \out(N)^{[ys]}\{2\}$.
\end{proof}

Let $\overline{\CI}\subset\CW/s$ be an upwardly closed subset. Then its preimage $\CI\subset\CW$ is upwardly closed as well and  Proposition \ref{coor} implies that $\on M\in\CV^\prime$ if $M\in\CV$. In order to prove that also $\out$ preserves Verma flags  we need a local definition of $\out$. 

\subsection{The local definition of translation out of the wall}

Let $\CG$ and $\CG/s$ be the moment graphs associated to $\Lambda$ and $\Lambda^\prime$, respectively.  Let $N\in\CV^\prime$. By Lemma \ref{WonZ} the $\CZ$-module $\out N$ is graded free over $S$ and, by \cite[Proposition 4.6]{Fie05}, it coincides with the global sections of the sheaf $\SM:=\Loc(\out N)$ on $\CG$. In this section we want to describe $\SM$. Let $\SN=\Loc(N)$ be the sheaf on $\CG/s$ associated to $N$.

\begin{lemma}\label{locout} For a vertex $x$ of $\CG$ we have $\SM^x=\SN^{\ol x}$. If $E\colon x\linie y$ is an edge and $y\neq xs$, let $\ol E\colon \ol x\linie\ol y$ be its image in $\CG/s$. Then we have $\SM^E=\SN^{\ol E}$ and $\rho_{x,E}\colon\SM^x\to\SM^E$ is identified with $\rho_{\ol x,\ol E}\colon \SN^{\ol x}\to\SN^{\ol E}$. If $E\colon x\stackrel{\alpha}\llinie xs$, then $\SM^E=\SN^{\ol x}/ \alpha\cdot \SN^{\ol x}$ and $\SM^x\to\SM^E$ equals the canonical quotient map $\SN^{\ol x}\to \SN^{\ol x}/ \alpha\cdot \SN^{\ol x}$.
\end{lemma}

\begin{proof} Choose $x\in\CW$. By definition, $\SM^x=(\out N)^{x}$, and by Proposition \ref{coor}, $(\out N)^{x}=N^{\ol x}=\SN^{\ol x}$. Let $E\colon x\linie y$ be an edge and suppose that $\ol x\neq\ol y$. Now $(\out N)^{x,y}$ is the image of $\out N$ in $(\out N)^x\oplus (\out N)^y=N^{\ol x}\oplus N^{\ol{y}}$. We want to show that it equals $N^{\ol x,\ol y}$, the image of  $N$ inside $N^{\ol x}\oplus N^{\ol{y}}$. According to Lemma \ref{WonZ}, we have a decomposition $\out N=N\oplus c^s\cdot N$. The image of the first summand, $N$, in $N^{\ol x}\oplus N^{\ol{y}}$ is $N^{\ol x,\ol{y}}$. We claim that the image of the second summand, $c^s\cdot N$, is contained in the image of the first. 

Suppose that $y=tx$. 
The image of an element in $c^s\cdot N$ is of the form $(x(\alpha_s)\cdot m,tx(\alpha_s)\cdot n)$ for some $(m,n)\in N^{\ol x,\ol{y}}$. Since $x(\alpha_s)\equiv tx(\alpha_s)\mod\alpha_t$ it is enough to show that $(\alpha_t\cdot m,0)\in  N^{\ol x,\ol{y}}$. By Lemma \ref{loccen} there is $z\in\CZ^{\prime}$ with $z_{\ol x}=\alpha_t$ and $z_{\ol{tx}}=0$, hence acting with $z$ on $(m,n)$ proves our claim. So we showed that $(\out N)^{x,y}=N^{\ol x,\ol y}$ and from the definition of $\Loc$ it follows that $\SM^E=\SN^{\ol E}$. 

Finally choose an edge $E\colon x\stackrel{\alpha}\llinie xs$. 
By Proposition \ref{coor}, $(\out N)^{x,xs}=\out(N^{\ol x})=\CZ^{x,xs}\otimes_S N^{\overline{x}}$. The latter space is canonically identified with $\{(m,n)\in N^{\ol x}\oplus N^{\ol x}\mid m-n\in \alpha\cdot N^{\ol x}\}$,  and we deduce $\SM^E= N^{\ol x}/\alpha\cdot N^{\ol x}= \SN^{\ol x}/\alpha\cdot\SN^{\ol x}$. 
\end{proof}

\begin{proposition} Let $N\in\CV^\prime$. Then $\out N\in\CV$. 
\end{proposition}
\begin{proof} According to Lemma \ref{loccar} it is enough to show that $\SM:=\Loc(\out N)$ is flabby and that the subquotients $\SM^{[x]}$ are graded free. 

Let $\Omega\subset\CW$ be an upwardly closed subset. We have to show that $\Gamma(\SM)\to\Gamma(\SM|_\Omega)$ is surjective. First suppose that $\Omega$ is $s$-invariant, i.e.~ $\Omega\cdot s=\Omega$. Then $\SM|_\Omega=\Loc((\out N)^{\Omega})=\Loc(\out(N^{\ol\Omega}))$.  Since $\out(N^{\ol\Omega})$ is free over $S$, $\Gamma(\SM|_\Omega)=(\out N)^{\Omega}$. The restriction $\Gamma(\SM)\to\Gamma(\SM|_\Omega)$ is hence identified with $\out N\to(\out N)^\Omega$ and, in particular, is surjective.

By \cite{Fie05}, Proposition 4.2, in order to show that $\SM$ is flabby it is enough to show that the map $\Gamma(\{\geq y\},\SM)\to\Gamma(\{> y\},\SM)$ is surjective for any $y\in\CW$. Let $m\in\Gamma(\{>y\},\SM)$. If $ys>y$, then the set $\{\geq y\}\setminus \{y,ys\}$ is $s$-invariant, hence the restriction of $m$ to $\{\geq y\}\setminus\{y,ys\}$ extends to a global section. So  we can assume that $m$ is supported on $\{ys\}$. Then we can extend $m$ to the vertex $y$ by setting $m_y:=m_{ys}$ (recall that $\SM^y=\SM^{ys}$ by construction).

If $ys<y$ we can assume, per induction, that $m$ extends to each vertex $z$ with $z>ys$, $z\neq y$ and we get a section $m^\prime\in\Gamma(\{\geq ys\}\setminus\{y,ys\},\SM)$. But $\{\geq ys\}\setminus\{y,ys\}$ is $s$-invariant, hence $m^\prime$ extends to a global section. In particular, $m$ extends to the vertex $y$.

So we showed that $\SM$ is flabby. Now $\SM^{[x]}=(\out N)^{[x]}$ and Proposition \ref{coor} shows that this space is graded free. 
\end{proof}

Hence we showed that $\on$ and $\out$ induce functors $\on\colon\CV\to\CV^\prime$ and $\out\colon\CV^\prime\to\CV$.

\begin{proposition} 
\begin{enumerate} 
\item The functors $\on\colon\CV\to\CV^\prime$ and $\out\colon\CV^\prime\to\CV$ are exact.
\item There are equivalences $\on\circ D\cong D\circ\on$ and $D\circ\out\cong\{2\}\circ\out\circ D$.
\end{enumerate}
\end{proposition}

\begin{proof}
The exactness of $\on$ and $\out$ follows from Proposition \ref{coor} and Lemma \ref{def exactness}. 
Since $\on M= M$ as an $S$-module it is clear that $\on$ commutes with the duality. For $N\in\CV^\prime$ we have $\out(DN)\{2\}=\CZ\{2\}\otimes_{\CZ^s} D N$ and,  using Proposition \ref{adj},
\begin{eqnarray*}
\CZ\{2\}\otimes_{\CZ^s} D N & = &  \Hom_{\CZ^s}(\CZ,\Hom_S(N,S)) \\
&  \cong & \Hom_{S}(\CZ\otimes_{\CZ^s}  N, S)\\
& = & D(\out N).
\end{eqnarray*}
\end{proof}

As a summarizing corollary of all the results in this section we get the following.
\begin{corollary}\label{ex and prop of trans} The functor $\vartheta_s:=\out\circ\on\circ\{1\}\colon \CV\to\CV$ has the following properties:
\begin{enumerate}
\item $\vartheta_s$ is exact and self-adjoint.
\item $\vartheta_s$ commutes with the duality.
\item For $\CI\subset\CW$ with $\CI=\CI\cdot s$ and $M\in\CV$ we have $(\vartheta_s M)^{\CI}=\vartheta_s(M^{\CI})$ and $(\vartheta_s M)_{\CI}=\vartheta_s(M_{\CI})$. 
\item For $y\in\CW$ with $ys<y$ we have 
$
(\vartheta_s M)^{[y]}\{-1\} \cong M^{[ys,y]} \cong (\vartheta_s M)^{[ys]}\{1\}.
$
\end{enumerate} 
\end{corollary}

\subsection{Translation of $S$-bimodules} Let $V$ be a reflection faithful representation, $\Lambda\subset V$ a regular orbit, and recall the equivalences  $F\colon \CV\stackrel{\sim}\to\CF_\Nabla$ and $G\colon \CF_\Nabla\stackrel{\sim}\to\CV$. For a simple reflection $s\in\CS$ Soergel defined a translation functor $\theta_s\colon\CF_\Nabla\to\CF_\Nabla$ by $\theta_s(M)=S\{1\}\otimes_{S^s}M$. 

\begin{proposition}\label{bimod trans} There are isomorphisms $G\circ\theta_s\cong \vartheta_s \circ G$ and $F\circ\vartheta_s\cong \theta_s \circ F$.
\end{proposition}
 \begin{proof}  Recall the map $\sigma\colon S\to \CZ$ defined in Section \ref{S-bimod}. It commutes with the action of $\CW$ on both algebras, hence induces a map $S^s\to\CZ^s$. By \cite[Lemma 4.5]{Soe04} we have $S=S^s\oplus\alpha_s\cdot S^s$. This decomposition is compatible with the decomposition $\CZ=\CZ^s\oplus c^s\cdot \CZ^s$, hence $S\{1\}\otimes_{S^s} F(M)\cong F(\CZ\{1\}\otimes_{\CZ^s} M)$ for $M\in\CV$, hence $ \theta_s \circ F\cong F\circ\vartheta_s$. The statement for $G$ follows. 
\end{proof}

\section{Projective objects}

Let $V$ be a reflection faithful representation of $(\CW,\CS)$ and let $\Lambda=\CW.v$ be a regular orbit. We identify $\Lambda$ with $\CW$ and provide $\Lambda$ with the Bruhat order. Let $\CZ=\CZ(\Lambda)$ and $\CV=\CV(\Lambda)$. Let $l\colon\CW\to\DN$ be the length function associated to $\CS$. For $M\in\CV$ define $\supp M:=\{w\in\CW\mid M^w\neq 0\}$. Recall that an object $P\in\CV$ is called projective, if the functor $\Hom_{\CV}(P,\cdot)$ maps short exact sequences in $\CV$ to short exact sequences of $k$-vector spaces. 

\begin{theorem}\label{B(x)}
For all $x\in\CW$ there exists an object $ B(x)\in\CV$, unique up to isomorphism, with the following properties:
\begin{enumerate}
\item \label{indec} $ B(x)$ is indecomposable and projective in $\CV$. 
\item \label{supp} $\supp  B(x)\subset\{y\mid y\leq x\}$ and $ B(x)^x\cong S\{l(x)\}$.
\end{enumerate}
Moreover $B(x)$ is self-dual, i.e.~$D(B(x))\cong B(x)$ as a $\CZ$-module. Each projective object in $\CV$ is isomorphic to a finite direct sum of modules of the form $B(x)\{k\}$ for $x\in\CW$ and $k\in\DZ$. 
\end{theorem}

\begin{proof} The following construction of $B(x)$ is well-known in similar situations (cf.~ \cite{Soe90,Fie03,Soe04}). We prove the existence of $ B(x)$ by induction on the Bruhat order. Obviously $B(e)=M(e)$ has all the desired properties. Let $x\in \CW$ and $s\in\CS$ such that $xs<x$ and suppose that we have already constructed $ B(xs)$. Consider $\vartheta_{s} B(xs)\in\CV$. By adjointness, $\Hom(\vartheta_{s} B(xs),\cdot)=\Hom( B(xs),\vartheta_{s}\cdot)$.
Since $\vartheta_s$ is exact and $ B(xs)$ is projective, $\vartheta_{s} B(xs)$ is projective.

From Corollary \ref{ex and prop of trans} we deduce $\supp \vartheta_s B(xs)\subset\{\leq xs\}\cup\{\leq xs\}\cdot s=\{\leq x\}$ and $(\vartheta_{s} B(xs))^x\cong B(xs)^{xs}\{1\}\cong S\{l(x)\}$.
Any 
indecomposable direct summand $ B(x)$ of $\vartheta_s B(xs)$ with $B(x)^x\cong S\{l(x)\}$ has properties (\ref{indec}) and (\ref{supp}). Hence we showed the existence of $ B(x)$. In order to show the uniqueness of $B(x)$ we need the following lemma.
\begin{lemma} 
If $f\in\End( B(x))$ is such that $f^x\in\End(B(x)^x)$ is an automorphism, then $f$ is an automorphism.
\end{lemma}
\begin{proof}
Applying the Fitting decomposition to any graded component of $ B(x)$ shows that 
$ B(x)=\varprojlim\im f^k\oplus\varinjlim\ker f^k$.
Since $ B(x)$ is indecomposable and since $f^x$ is an automorphism, $\ker f=0$.
\end{proof}

Suppose that $ B(x)^\prime$ is another object having properties (\ref{indec}) and (\ref{supp}). By projectivity we can lift the isomorphisms
$( B(x)^\prime)^x\cong S\{l(x)\}\cong ( B(x))^x$
to maps $f\colon B(x)\to B(x)^\prime$ and $g\colon B(x)^\prime\to B(x)$. The composition $g\circ f\colon  B(x)\to  B(x)$ is such that $(g\circ f)^x\colon B(x)^x\to B(x)^x$ is the identity. Hence $g\circ f$ is an automorphism by the lemma. So $B(x)$ is a direct summand of $ B(x)^\prime$, hence $B(x)\cong B(x)^\prime$. A similar argument shows that every projective object in $\CV$ is isomorphic to a direct sum of shifted copies of various $B(x)$'s.  

Suppose we have already shown that the object $B(y)$ is self-dual for all $y<x$. The object $\vartheta_{s} B(xs)$ is  self-dual by Corollary \ref{ex and prop of trans} and $\vartheta_s B(xs)\cong B(x)\oplus R$, where $R$ is isomorphic to a direct sum of modules of the form $B(y)\{k\}$ with $y<x$ and $k\in\DZ$. Choose an isomorphism
$D(\vartheta_s B(xs))\cong \vartheta_s B(xs)$.
This induces a map $D( B(x))\to B(x)$ on the direct summands. Now $R^x=0$ and hence $(DR)^x=0$, so $D( B(x))^x\to B(x)^x$ is an isomorphism.

$D(B(x))$ admits a Verma flag since it is a direct summand of $D(\vartheta_s B(xs))$. By projectivity of $B(x)$ we deduce, as before, that $B(x)$ is a direct summand of $D(B(x))$. Since the latter is indecomposable we get $D( B(x))\cong B(x)$, which was left to be shown.
\end{proof}

Let $\CP\subset\CV$ be the full subcategory of all projective objects. Soergel defined a full subcategory $\CB\subset\CF_\Nabla$ of ``special'' objects with prescribed classes in the Grothendieck group of $\CF_\Nabla$. He then showed that the objects in $\CB$ are the direct sums of direct summands of objects of the form $\vartheta_{s}^\prime\cdots\vartheta_t^\prime( S(e)\{k\})$ for a sequence $(s,\dots,t)$ of simple reflections and $k\in\DZ$. The proof of the theorem above together with Proposition \ref{bimod trans} shows the following.  
\begin{theorem}\label{BequivP} The equivalence $\CF_\Nabla\cong\CV$ induces an equivalence $\CB\cong\CP$. 
\end{theorem}

\subsection{The Braden-MacPherson construction}\label{BMP construction}
We give a second construction of $B(x)$ as the global sections of the ``canonical sheaf'' of Braden and MacPherson (cf.~ \cite{BMP}). Let $\CG=\CG(\Lambda)$ and denote by $\SB(x)$ the canonical sheaf associated to the full subgraph $\CG_{\leq x}\subset\CG$ whose set of vertices is $\{y\mid y\leq x\}$. 

We recall the construction of $\SB(x)=\{\SB(x)^y,\SB(x)^E,\rho_{y,E}\}$. We set $\SB(x)^x=S$ and $\SB(x)^y=0$ for all $y\not\leq x$. Suppose we already have constructed $\SB(x)^y$ for some $y$, and let $E\colon y^\prime\stackrel{\alpha}\to y$ be an edge (in particular, $y^\prime< y$). Set $\SB(x)^E=\SB(x)^y/ \alpha\cdot\SB(x)^y$ and let $\rho_{y, E}\colon \SB(x)^y\to \SB(x)^E$ be the canonical map.  

Let $y\in\CW$ and suppose that we have constructed $\SB(x)$ on the subgraph $\CG_{>y}$. We have to construct $\SB(x)^y$ together with all maps $\rho_{y,E}\colon\SB(x)^y\to\SB(x)^E$ for edges  $E\in\CE^{\delta y}:=\{E\mid E\colon y\to w\}\subset\CE$. Define the sections of $\SB(x)$ on $\CG_{>y}$ as
$$
\Gamma(\{>y\},\SB(x)):=\left\{(m_{w})\in\bigoplus_{w>y} \SB(x)^{w}\left| \,
\begin{matrix}
\rho_{w,E}(m_w)=\rho_{w^\prime,E}(m_{w^\prime}) \text{ for any} \\
\text{edge $E\colon w\linie w^\prime$ of $\CG_{>y}$}
\end{matrix}
\right.
\right\}
$$
and let  $\SB(x)^{\delta y}\subset \bigoplus_{E\in\CE^{\delta y}}\SB(x)^E$ be the image of 
$$
\bigoplus{\rho_{w,E}}\colon \Gamma(\{>y\},\SB(x))\to\bigoplus_{E\in\CE^{\delta y}}\SB(x)^E.
$$
Let $\SB(x)^y\to\SB(x)^{\delta y}$ be a projective cover in the category of graded $S$-modules, and let $\rho_{y,E}\colon\SB(x)^y\to\SB(x)^E$ be its components. This finishes the construction of $\SB(x)$. 

Recall that we defined $\SB(x)^{[y]}$ as the kernel of the map $\SB(x)^y\to\SB(x)^{\delta y}$. From the construction alone it is not clear wether $\Gamma(\SB(x))\in\CV$, i.e.~ whether $\SB(x)^{[y]}$ is graded free for all $y$.  But one can easily prove a weaker statement. Let $\CZ\catmod^{ref}\subset\CZ\catmod^f$ be the full subcategory of objects $M$ such that for any upwardly closed $\CI\subset\CW$ the quotient $M^{\CI}$ is {\em reflexive} as an $S$-module, i.e.~ isomorphic to its double dual. There is an analogously defined exact structure on $\CZ\catmod^{ref}$.
\begin{theorem}[{\cite[Theorem 5.2]{Fie05}}] $\Gamma(\SB(x))$ is an object in $\CZ\catmod^{ref}$. It is characterized, up to isomorphism, by the following properties:
\begin{enumerate}
\item $\Gamma(\SB(x))$ is indecomposable and projective in $\CZ\catmod^{ref}$.
\item $\supp\Gamma(\SB(x))\subset\{y\mid y\leq x\}$ and $\Gamma(\SB(x))^x\cong S$.
\end{enumerate}
\end{theorem}

\begin{corollary}\label{locglob} We have $\Gamma(\SB(x))\cong B(x)\{-l(x)\}$. In particular, $\Gamma(\SB(x))\in\CV$. 
\end{corollary}
\begin{proof}  Proposition \ref{coor} shows that $\vartheta_s$ preserves $\CZ\catmod^{ref}$ as well, and the construction of $B(x)$ shows that $B(x)$ is projective even in $\CZ\catmod^{ref}$. 
\end{proof}

For $y\in\CW$ there are induced isomorphisms $\SB(x)^{[y]}\cong B(x)^{[y]}\{-l(x)\}$ and $\SB(x)^y\cong B(x)^y\{-l(x)\}$.

\section{Properties of $B(x)$}
Choose $x,y\in\CW$ and consider the canonical inclusions $B(x)_y\subset B(x)^{[y]}\subset B(x)^y$. 

\begin{proposition}\label{struc of proj} Suppose that $y\leq x$.
\begin{enumerate}
\item\label{proj free} The $S$-modules $B(x)_y$, $B(x)^{[y]}$ and $B(x)^y$ are graded free over $S$. 
\item\label{proj homo} Let $\alpha_1,\dots,\alpha_{l(y)}$ be the labels of the edges ending at $y$ (note that $l(y)=\{t\in\CT\mid ty<y\}$). Then $B(x)_y=\alpha_1\cdots \alpha_{l(y)}\cdot B(x)^{[y]}$.
\item\label{proj dual} There is an isomorphism $B(x)^y\cong D(B(x)_y)$ of graded $S$-modules.
\item\label{proj deg} 
 Let $B(x)^y\cong \bigoplus_i S\{k_i\}$. 
Then $B(x)^{[y]}\cong\bigoplus_i S\{2l(y)-k_i\}$.
\item\label{proj cov} $B(x)^y\to B(x)^y/B(x)^{[y]}$ is a projective cover.
\end{enumerate}
\end{proposition}
\begin{proof} 
We already know that $B(x)\in\CV$, hence, by definition, $B(x)^{[y]}$ is graded free. By the Braden--MacPherson construction, $B(x)^y$ is graded free and (\ref{proj dual}) implies that $B(x)_y$ is graded free.

Recall that $\CE^{\delta y}$ is  the set of edges that originate at $y$. Let $\CE^y$ be the set of all edges linked to $y$. We have
$$
B(x)^{[y]}  =  \ker\left(\SB(x)^y\to\bigoplus_{E\in\CE^{\delta y}}\SB(x)^{E}\right)
$$ 
and
$$
B(x)_{y}  =  \ker \left(\SB(x)^y\to\bigoplus_{E\in\CE^y}\SB(x)^{E}\right).
$$ 
The set $\CE^y\setminus\CE^{\delta y}$ is the set of edges ending at $y$. For any such $E$ we have $\SB(x)^E=\SB(x)^y/\alpha\cdot \SB(x)^y$, where $\alpha$ is the label of $E$. Since no edge starting at $y$ is labelled by a multiple of one of the $\alpha_i$, every basis element of $B(x)^{[y]}$ (as a graded free $S$-module) is not divisible by $\alpha$ in $B(x)^y$, hence maps to a non-zero element in $\SB(x)^E$. So 
$B(x)_y=\alpha_1\dotsm \alpha_{l(y)}\cdot B(x)^{[y]}$, as claimed in (\ref{proj homo}).

Note that  $\Hom_\CV(B(x), M(y)\{k\})=\Hom_S(B(x)^y,S\{k\})$, hence $D(B(x)^y)$ is isomorphic to $\bigoplus_{k\in\DZ}\Hom_\CV(B(x), M(y)\{k\})$.  
Analogously, $B(x)_y$ is isomorphic to $ \bigoplus_{k\in\DZ}\Hom_\CV(M(y)\{-k\}, B(x))$. Since $D(M(y)\{k\})\cong M(y)\{-k\}$ there is an induced isomorphism $\Hom_\CV(M(y)\{-k\}, B(x))\cong\Hom_\CV(DB(x),M(y)\{k\})$ for any $k$, and the latter space is isomorphic to $\Hom_\CV(B(x),M(y)\{k\})$, hence (\ref{proj dual}).

Claim (\ref{proj deg}) follows immediately from (\ref{proj homo}) and (\ref{proj dual}). Claim (\ref{proj cov}) follows from Corollary \ref{locglob} and the construction of $\SB(x)$. 
\end{proof}

One of the main problems in Kazhdan--Lusztig theory is the following conjecture.
\begin{conjecture}\label{bound conjecture}
Let $x,y\in\CW$ with $y<x$. Then the $S$-module $B(x)^{[y]}$ lives in degrees $> -l(y)$. 
\end{conjecture}
\begin{remarks}
\begin{enumerate}
\item  The conjecture is proven for crystallographic Coxeter systems and fields $k$ of characteristic zero using the fact that  $B(x)$ is isomorphic to the intersection cohomology of a Schubert variety. The conjecture translates into one of the defining axioms for intersection cohomology.
\item  In the case that the characteristic of the ground field is positive (and big enough such that a reflection faithful representation exists) there is no direct proof of the conjecture. However, it follows from its characteristic zero analog if the characteristic is ``big enough'', though no explicit bounds are known so far. Very little is known about the conjecture in small characteristics.  
\item  One finds a proof of the conjecture for the dihedral cases in \cite{Soe04}. We will give a proof for universal Coxeter systems in the last section of this paper. Both results hold in arbitrary characteristic. 
\end{enumerate}
\end{remarks}

\begin{proposition}\label{local struc of proj} Let $x,y\in\CW$ and $s\in\CS$ such that $ys<y$. Then $B(x)^{[y,ys]}$ is isomorphic to a direct sum of shifted copies of $M(ys)$ and $P(ys,y)$. 
\end{proposition}
\begin{proof} Note that $B(x)^{[y,ys]}$ is graded free over $S$ since it is an extension of $B(x)^{[y]}$ and $B(x)^{[ys]}$. By Lemma \ref{local structure} we only have to show that no shifted copy of $M(y)$ occurs. Let $E\colon ys\stackrel{\alpha}\to y$ be the edge connecting $y$ and $ys$. By the Braden-MacPherson construction we have $\SB(x)^E=\SB(x)^{y}/ \alpha\cdot \SB(x)^{y}$, hence any copy of $M(y)$ has to extend with some $M(ys)$, as no generator of $\SB(x)^{[y]}$ is divisible by $\alpha$ in $\SB(x)^y$, so does not map to zero in $\SB(x)^{E}$.
\end{proof}

For a graded $S$-module $N$ we denote by  $N_{\{\leq m\}}\subset N$ the submodule generated in degrees $\leq m$, i.e.~ $N_{\{\leq m\}}=\sum_{j\leq m}S\cdot N_{\{j\}}$.

\begin{proposition} \label{purity - bound}
\begin{enumerate}
\item\label{bound} Let $x\in\CW$ and $s\in\CS$ such that $xs<x$. Assume that Conjecture \ref{bound conjecture} holds for $xs$, i.e.~ that $B(xs)^{[y]}$ lives in degrees $> -l(y)$ for all $y<xs$. Then $B(x)^{[y]}$ lives in degrees $\geq -l(y)$  for all $y<x$. 
\item \label{bound - ind} Choose $x,y\in\CW$ such that $y<x$ and assume that $B(x)^{[y]}$ lives in degrees $\geq -l(y)$. Then $B(x)^{y}$ is generated in degrees $\leq -l(y)$. Suppose furthermore that 
$$
 B(x)^{[y]}_{\{\leq -l(y)\}}\cap B(x)^y_{\{\leq -l(y)-1\}}=0.
$$
Then $B(x)^{[y]}_{\{\leq -l(y)\}}=0$, i.e.~ Conjecture \ref{bound conjecture} holds for the pair $x,y$.
\end{enumerate}
\end{proposition}

\begin{proof} In order to show (\ref{bound}) it suffices to show that $(\vartheta_s B(xs))^{[y]}$ lives in degrees $\geq -l(y)$ since $B(x)^{[y]}$ is a direct summand.  Assume that $ys<y$. Then $B(xs)^{[y,ys]}$ lives in degrees $>-l(y)$. By Corollary \ref{ex and prop of trans} we have $(\vartheta_s B(xs))^{[y]}\cong B(xs)^{[y,ys]}\{1\}$, hence $(\vartheta_s B(xs))^{[y]}$ lives in degrees $\geq -l(y)$, and $(\vartheta_s B(xs))^{[ys]}\cong B(xs)^{[y,ys]}\{-1\}$ lives in degrees $>-l(y)+1=-l(ys)$.

Now we will prove (\ref{bound - ind}). The first claim follows from Proposition \ref{struc of proj}. The map $f\colon B(x)^{[y]}_{\{\leq -l(y)\}}   \to  B(x)^y/  B(x)^y_{\{\leq -l(y)-1\}}$ is injective by the second assumption. From Proposition \ref{struc of proj} and our first assumption we deduce that both modules are isomorphic to a direct sum of copies of $S\{l(y)\}$ and have the same rank. Hence $f$ is an isomorphism. Now $ B(x)^{[y]}_{\{\leq -l(y)\}}$ being non-zero contradicts  the fact that  $B(x)^y\to B(x)^{ y}/B(x)^{[y]}$ is a projective cover. 
\end{proof}

Choose $x,y\in \CW$ with $y<x$ and let $\{k_i\}_i$ be the multiset of numbers such that $B(x)^{[y]}\cong\bigoplus_i S\{k_i\}$. Conjecture \ref{bound conjecture} claims that $k_i<l(y)$ for all $i$. We have $B(x)^y\cong\bigoplus_iS\{l_i\}$ with $l_i=2l(y)-k_i$. 
From Corollary \ref{locglob} and the inductive construction of $\SB(x)$ it follows that one of the $l_i$ equals $l(x)$ and the others are strictly smaller than $l(x)$. Hence Conjecture \ref{bound conjecture} is equivalent to $2l(y)-l(x)\leq k_i <l(y)$ for all $i$. By induction and Proposition \ref{purity - bound} we can assume that $2l(y)-l(x)\leq k_i \leq l(y)$.  The following lemma rules out a very specific case and is used to prove the conjecture for universal Coxeter systems in Section \ref{univCox}.

\begin{lemma}\label{simple case}  $B(x)^{[y]}$ is never isomorphic to $S\{l(y)\}\oplus S\{2l(y)-l(x)\}$. 
\end{lemma}

\begin{proof} Suppose that $B(x)^{[y]}$ is isomorphic to $S\{l(y)\}\oplus S\{2l(y)-l(x)\}$. Then $\SB(x)^{[y]}\cong S\{l(y)-l(x)\}\oplus S\{2l(y)-2l(x)\}$ and  $\SB(x)^y=S\{l(y)-l(x)\}\oplus S$. The inclusion $\SB(x)^{[y]}\to \SB(x)^y$ gets identified with a map 
$$
S\{l(y)-l(x)\}\oplus S\{2l(y)-2l(x)\}\to  S\{l(y)-l(x)\}\oplus S.
$$ 

Let $(f,g)\in S\{l(y)-l(x)\}\oplus S$ be the image of a generator of the copy of $S\{l(y)-l(x)\}$ on the left hand side. If $f\neq 0$, then the induced map $S\{l(y)-l(x)\}\to S\{l(y)-l(x)\}$ is an isomorphism. But this implies that there is a direct summand in $\SB(x)^y$ lying in the kernel of $\SB(x)^y\to\SB(x)^{\delta y}$, which contradicts the fact that this mapping is a projective cover. Hence $f=0$. We deduce that there is an element $g\in S$ of degree $l(x)-l(y)$ such that $g\cdot m_y\in\SB(x)^{[y]}$, where $m_y\in\SB(x)^y$ is a non-zero element of degree zero. 

Consider a non-zero element $m_x$ in $\SB(x)^x\cong S$ of degree zero. By construction of $\SB(x)$ it extends to a global section $(m_v)\in \Gamma(\SB(x))$ such that  each $m_v\in\SB(x)^v$ is a generator of the degree zero component. In particular, for any edge $E\colon y\to y^\prime$ we have $\rho_{y,E}(m_y)\neq 0$ in $\SB(x)^E$. So the degree of the map $g$ constructed above is at least twice the number of edges $E\colon y\to y^\prime$ with $y^\prime\leq x$. By Deodhar's conjecture, proven in \cite{Dyer1}, this number is at least $l(x)-l(y)$, so $\deg g\geq 2(l(x)-l(y))$, and we have a contradiction.   
\end{proof}

\section{Combinatorics in the Hecke algebra}

Let $\CH=\CH(\CW,\CS)$ be the Hecke algebra of $(\CW,\CS)$, i.e.~   the free $\DZ[v,v^{-1}]$-module with basis $\{T_x\}_{x\in\CW}$ that is endowed with a multiplication such that \begin{eqnarray*}
{ T}_x\cdot { T}_{y} & = & { T}_{xy}\quad\text{if $l(xy)=l(x)+l(y)$}, \\
{ T}_s^2 & = & v^{-2}T_e+(v^{-2}-1)T_s \quad\text{for $s\in\CS$}.
\end{eqnarray*}
${ T}_e$ is a unit in $\CH$ and for any $x\in\CW$ there exists an inverse of
${ T}_x$ in $\CH$. For $s\in\CS$ we have ${ T}_s^{-1}=v^2T_s+(v^2-1)$. There is a duality (i.e.~ a $\DZ$-linear anti-involution) $d\colon\CH\to\CH$, given by 
$d(v)=v^{-1}$ and $d(T_x)  =  T_{x^{-1}}^{-1}$ for $x\in\CW$.

\subsection{The self-dual Kazhdan-Lusztig basis}
 Set $\tilde T_x=v^{l(x)}  T_x$.

\begin{theorem}[\cite{KL,Soe97}]\label{self-dual elts} For any $x\in\CW$ there exists a unique element $C_x^\prime=\sum_{y\in\CW} h_{y,x}(v)\cdot { \tilde T}_y\in\CH$ with the following properties:
\begin{enumerate}
\item \label{prop of Cs: duality} $C_x^\prime$ is self-dual, i.e.~ $d(C_x^\prime)=C_x^\prime$. 
\item \label{prop of Cs: support} $h_{y,x}(v)=0$ if $y\not\leq x$, and $h_{x,x}(v)=1$,
\item \label{prop of Cs: norm} $h_{y,x}(v)\in v\DZ[v]$ for $y<x$.
\end{enumerate}
\end{theorem} 

For each pair $x,y$ there is a unique polynomial $P\in\DZ[v,v^{-1}]$ such that $P_{y,x}(v^{-2})=v^{l(y)-l(x)}h_{y,x}(v)$. These are the well-known {\em Kazhdan-Lusztig polynomials}.

 We want to recall the main idea of the inductive construction of $C^\prime_x$. We begin with $C^\prime_e=\tilde T_e$ and $C^\prime_s= \tilde T_s+v$ for $s\in\CS$. Let $x\in\CW$, $s\in\CS$ with $xs<x$ and assume that $C^\prime_{xs}$ is already constructed. Define $b_{y}\in\DZ[v,v^{-1}]$ by $C^\prime_{xs}\cdot C^\prime_s=\sum_y b_y(v)\cdot \tilde T_y$. Then $b_y\in\DZ[v]$ for all $y\in\CW$ and 
$$
C^\prime_x=C^\prime_{xs}\cdot C^\prime_s-\sum_{y<x} b_y(0)\cdot C^\prime_y.
$$

\begin{proposition}\label{calc in Hecke} Let $x,y\in\CW$, $s\in\CS$ with $xs<x$ and $y<x$. 
\begin{enumerate}
\item\label{bs} If $ys<y$, then $v\cdot b_y(v)=b_{ys}(v)=v\cdot h_{ys,xs}(v)+h_{y,xs}(v)$.
\item\label{subcritical elts} If $b_y(0)\neq 0$, then $ys<y$, and for any $t\in\CS$ with $xst<xs$ either $yt<y$ or $y=xst$ holds.
\item\label{prop of KL polys} If  $ys<y$, then $h_{ys,x}(v)=v\cdot h_{y,x}(v)$.\end{enumerate}
\end{proposition}

\begin{proof}
Let $y\in\CW$ and suppose  that $ys<y$. We claim that   
\begin{equation*}
v\cdot { \tilde T}_{y}\cdot C^\prime_s = { \tilde T}_{ys}\cdot C^\prime_s={ \tilde T}_y +  v\cdot { \tilde T}_{ys}.\eqno{(\ast)}
\end{equation*}
Recall that $C^\prime_s={ \tilde T}_s+v$.
The second identity of $(\ast)$  follows directly from the definition of the multiplication. Also, by definition, $(C_s^\prime)^2=(v^{-1}+v)C_s^\prime$. Hence
$$
{\tilde T}_{ys}\cdot (C_s^\prime)^2={\tilde T}_{ys}\cdot(v^{-1}+v)\cdot C_s^\prime
$$
and
$$
{\tilde T}_{ys}\cdot (C_s^\prime)^2=({\tilde T}_{y}+v\cdot {\tilde T_{ys}})\cdot C_s^\prime
$$ 
and from this the first identity follows. We deduce part (\ref{bs}) of the proposition directly from $(\ast)$.

We prove the first part of (\ref{subcritical elts}). Suppose $b_y(0)\neq0$. If $y<ys$, then, by part (\ref{bs}), $b_y(v)=v\cdot h_{y,xs}(v)+h_{ys,xs}(v)$, hence $h_{ys,xs}(0)\neq0$, hence $xs=ys$, which contradicts $y<x$. Hence $ys<y$.

Let us prove part (\ref{prop of KL polys}). Suppose $ys<y$. Then $v\cdot b_y(v)=b_{ys}(v)$ by  part (\ref{bs}). For all $z$ we get $h_{z,x}$ from $b_z$ by substracting $b_{x^\prime}(0)\cdot h_{z,x^\prime}$ for all $x^\prime<x$. Now $b_{x^\prime}(0)\neq 0$ implies $x^\prime s<x^\prime$ by the already proven part of (\ref{subcritical elts}). By induction (on the Bruhat order of $x$) we assume that $v\cdot h_{y,x^\prime}(v)=h_{ys,x^\prime}(v)$, hence the claim.  

Finally we prove the second part of (\ref{subcritical elts}). Assume  $b_y(0)\neq 0$. We already showed that $ys<y$. Then  $b_y(v)=h_{ys,xs}(v)+v^{-1}\cdot h_{y,xs}(v)$ by part (\ref{bs}) and we deduce $(v^{-1}\cdot h_{y,xs}(v))(0)\neq 0$. Suppose that $y<yt$. Then, by part (\ref{prop of KL polys}), applied to $xst<xs$, $h_{y,xs}(v)=v\cdot h_{yt,xs}(v)$. Hence $h_{yt,xs}(0)\neq 0$ and this implies $yt=xs$. 
\end{proof}

\subsection{Characters in $\CH$}
Let $N\cong\bigoplus_i S\{k_i\} $ be a graded free $S$-module. The multiset of numbers $\{k_i\}$ is well-defined and hence we can define the {\em graded character of $N$} by
$\udim(N)=\sum_i v^{-k_i}\in\DZ[v,v^{-1}]$.
For $ M\in\CV$ we define its {\em graded character} as 
$$
h( M)=\sum_{y\in\CW}\udim(M^{[y]})\cdot v^{l(y)}{ \tilde T}_y\in\CH.
$$
This defines a group homomorphism from the split Grothendieck group of $\CV$ to $\CH$.

\begin{conjecture}[cf.~ {\cite[Vermutung 1.13]{Soe04}}]\label{KL conjecture}
For any $x\in\CW$ we have
$h( B(x))=C_x^\prime$.
\end{conjecture}

For $x,y\in\CW$ set $f_{y,x}(v)=\udim(B(x)^{[y]})\cdot v^{l(y)}\in\DZ[v,v^{-1}]$. Then conjecture \ref{KL conjecture} amounts to showing that $f_{y,x}=h_{y,x}$. 

\begin{proposition}\label{eqcon} Conjecture \ref{bound conjecture} is equivalent to Conjecture \ref{KL conjecture}.
\end{proposition}
\begin{proof} Let $x\in\CW$. We consider the defining properties of $C^\prime_x$ in Theorem \ref{self-dual elts}. In any case $h(B(x))$ is self-dual, as follows from the self-duality of $B(x)$ and the  following lemma. Moreover, the support property follows from the support property of $ B(x)$ in Theorem \ref{B(x)}. And the normalization (\ref{prop of Cs: norm}) in Theorem \ref{self-dual elts} is a reformulation of Conjecture \ref{bound conjecture}. 
\end{proof}

For the proof of the following lemma one can copy Soergel's arguments or use the fact that the equivalence $\CB\cong\CP$ of Theorem \ref{BequivP} commutes with the  translation functors (cf.~ Proposition \ref{bimod trans}) and the duality. 
\begin{lemma}[\cite{Soe04}, Proposition 5.9 \& Bemerkung 6.16] Choose $s\in\CS$ and let $B\in\CV$ be projective. Then 
 $h(\vartheta_s B)=h( B)\cdot C_s^\prime$ and $h(D( B))=d(h( B))$.
\end{lemma}

\section{Universal Coxeter systems}\label{univCox}

Let $(\CW,\CS)$ be a universal Coxeter system, i.e.~ a Coxeter system such that for any $s,t\in\CS$ with $s\neq t$ the product $st$ has infinite order. In other words, $\CW$ is the group generated by $s\in\CS$ with the only relations $s^2=1$. The following lemma is an immediate consequence.

\begin{lemma}\label{red expr} Let $(\CW,\CS)$ be a universal Coxeter system. For any $w\in\CW$ there is a unique reduced expression
$w=s_1\cdots s_n$.  Hence for $w\neq e$ there is a unique simple root $s\in\CS$ such that $ws<w$.
\end{lemma} 

\begin{lemma}\label{critical elts for univ}
Choose $x\in\CW$ and $s\in\CS$ such that $xs<x$ and let $C^\prime_{xs}\cdot C^\prime_s=\sum b_y(v)  \cdot \tilde T_y$. Suppose there is $y<x$ with $b_y(0)\neq 0$. Then there is $t\in\CS$ such that $y=xst$ and $xsts<xst<xs<x$.
\end{lemma}

\begin{proof} We can assume that $x\neq s$. Then there exists $t\in\CS$  such that $xst<xs$. By Proposition \ref{calc in Hecke}, (\ref{subcritical elts})  we have $ys<y$, hence $y<yt$ by Lemma \ref{red expr}. Again using Proposition \ref{calc in Hecke}, (\ref{subcritical elts}) gives  $y=xst$.
\end{proof}

\begin{theorem} Suppose $(\CW,\CS)$ is universal. Then Conjecture \ref{KL conjecture} holds, i.e.~
$$
h( B(x))=C^\prime_x
$$
for any $x\in\CW$. 
\end{theorem} 
For the proof we need the following lemma which holds for arbitrary Coxeter systems.
\begin{lemma}\label{calculation} Choose $x\in\CW$ and $s,t\in\CS$ such that $xst<xs<x$. Then $h_{xs,x}(v)=v$ and  $h_{xst,x}(v)=v^2$. 
\end{lemma}
\begin{proof}  Let $C^\prime_{xs}\cdot C^\prime_s=\sum b_z(v)\cdot\tilde T_z$. Since $h_{xs,xs}(v)=1$ we have $b_{xs}(v)=v=h_{xs,x}(v)$, hence the first statement. We prove the second statement by induction on the Bruhat order of $x$. A simple calculation shows that it is true for $x=ts$. So suppose that the lemma is proven for all $x^\prime$ with $x^\prime<x$ (and arbitrary $s,t\in\CS$). First suppose that $xsts<xst$. Then, by the induction hypothesis, we have $h_{xst,xs}(v)=v$ and $h_{xsts,xs}(v)=v^2$. Then $b_{xst}(v)=1+v^2$ and hence $h_{xst,x}(v)=v^2$. If $xsts>xst$, then $xsts\not\leq xs$, hence $h_{xst,xs}(v)=v$ and $h_{xsts,xs}(v)=0$. Hence $b_{xst}(v)=v^2=h_{xst,x}(v)$. 
\end{proof}

\begin{proof}[Proof of the Theorem] 
We prove the theorem by induction on the Bruhat order. The case $x=e$ is clear. Suppose the theorem is proven for all $x^\prime < x$, i.e.~ $f_{y,x^\prime}=h_{y,x^\prime}$ for all $y\in\CW$.
 
As we have seen in the proof of Proposition \ref{eqcon} we have to show that $f_{y,x}\in v\DZ[v]$ for all $y<x$. 
Suppose that $y<x$ is such that this is not the case. Choose $s\in\CS$ such that $xs<x$ and consider $C^\prime_{xs}\cdot C_s^\prime=\sum_y b_y(v) \cdot \tilde T_y$. Then $b_y(0)\neq 0$ and by Lemma \ref{critical elts for univ} there is $t\in\CS$ such that $y=xst$ and $ys<y$. The proof of Lemma \ref{calculation} shows that $b_y(v)=1+v^2$. But then $f_{y,x}(v)=1+v^2$, which contradicts Lemma \ref{simple case}. 
\end{proof}

\end{document}